\documentclass[a4paper,12pt]{article}
\usepackage{amsmath, amssymb, amscd, amsxtra, amsthm}

  \newtheorem{theorem}{Theorem}
  \newtheorem{lemma}[theorem]{Lemma}
  \newtheorem{proposition}[theorem]{Proposition}
  \newtheorem{corollary}[theorem]{Corollary}
  \newtheorem{definition}[theorem]{Definition}

\theoremstyle{definition}
  \newtheorem{remark}[theorem]{Remark}

\newcommand{\epsi}{{\epsilon}}

\newcommand{\F}{{\cal F}}

\newcommand{\Z}{{\bf Z}}

\newcommand{\R}{{\bf R}}
\newcommand{\RR}{{\bf R}}
\newcommand{\C}{{\bf C}}
\newcommand{\D}{{\bf D}}

\renewcommand{\F}{{\bf F}}

\newcommand{\End}{{\mathrm{End}}}

\newcommand{\HH}{{\bf{H}}}

\newcommand{\pt}{\mathrm{pt}}
\newcommand{\ch}{\mathrm {ch}}
\newcommand{\sign}{{\mathrm{sign}}\,}
\newcommand{\caF}{{| S |}}
\newcommand{\caFi}{{| S_i |}}
\newcommand{\caS}{{| S |}}
\newcommand{\tT}{{\tilde T}}
\newcommand{\hT}{{T}}
\newcommand{\rc}{{s}}
\newcommand{\Groupp}{{G}}
\newcommand{\Irr}{{Ir}}

\newcommand{\Ker}{\mbox{ \rm Ker}\,}

\newcommand{\Spin}{\mbox{\rm Spin}\,}
\newcommand{\rank}{\mbox{\rm rank}\,}
\newcommand{\Tr}{tr}
\newcommand{\Hom}{\mbox{\rm Hom}}
\newcommand{\Map}{\mbox{\rm Map}}

\newcommand{\WH}{{\widehat \HH}}
\newcommand{\Ind}{\mbox{\rm Ind}\,}

\newcommand{\id}{\mbox{\rm id}\,}
\newcommand{\even}{\mathrm{even}}
\newcommand{\odd}{\mathrm{odd}}

\newcommand{\dpin}{{\Gamma}}
\newcommand{\pin}{\mathrm{Pin^-(2)}}

\newcommand{\bysame}%
 {\leavevmode\hbox to 3em{\hrulefill}\,}

\title{Equivariant maps between sphere bundles over tori
and $KO^*$-degree}
\author{M.~Furuta and Y.~Kametani}
\begin{document}
\maketitle
\abstract{We show a non-existence result for
some class of equivariant maps between
sphere bundles over tori.
The notion of equivariant $KO^*$-degree is used in the proof.
As an application to Seiberg-Witten theory we have a new inequality
$b_2^{+}(X)\ge -\sign(X)/8 + c(X)+\varepsilon(X)$ 
for a connected closed oriented spin $4$-manifold $X$
with indefinite intersection form,
where
$c(X)$ is a non-negative integer determined by the quadruple
cup product on $H^1(X ;\Z)$ and 
some maps in $KO$-theory induced from the Albanese map of $X$,
and where $1 \le \varepsilon(X) \le 3$.

\section{Introduction}
Let $\tilde{\R}$ be the $1$-dimensional real line with
the linear involution $-1$ and let $\tilde{\R}^p=\tilde \R^{\oplus p}$.
The involution on $\tilde \R$ then induces an involution on
the $n$-dimensional torus $\tT^n= (\tilde \R/\Z)^n$,
so it is a Real space \cite{Atiyah-R}.
Fix an non-negative integer $l$.
We denote by $V_0$, $W_0$ the product vector bundles
$\tT^n \times \tilde \R^{x}$, $\tT^n \times \tilde \R^{x+l}$
respectively.
The diagonal actions on $V_0$, $W_0$ give lifts of the
involution.

Let $V_1$ and $W_1$ be symplectic vector bundles
over the Real space $\tT^n$ in the sense of Dupont \cite{Dupont};
each of them is defined to be a complex vector bundle
together with an anti-linear lift 
of the involution such that its square is  $-1$ on each fiber.
Since the symplectic action 
defines a quaternionic structure on the fiber over a fixed point,
the rank of symplectic vector bundles over $\tT^n$ should be even.
We thus denote $\rank_{\C}V_1-\rank_{\C}W_1=2k$.

Let $\pin$ be the subgroup of $Sp(1)$ generated
by $U(1)$ and $j$.
Then the quotient $\pin/U(1)$ is isomorphic to $\{ \pm 1\}$,
so we regard $V_0$, $W_0$ as $\pin$-equivariant bundles.
We define $\pin$-actions on $V_1$, $W_1$ 
by using complex scalar multiplication and the symplectic action.
Then our purpose of this paper is to give a necessary condition 
for the existence of 
a $\pin$-equivariant fiber-preserving map from 
the sphere bundle $S(V_0 \oplus V_1)$ to 
the sphere bundle $S(W_0 \oplus W_1)$,
or equivalently the existence of
a proper $\pin$-equivariant fiber-preserving map from 
$V_0 \oplus V_1$ to $W_0 \oplus W_1$,
under some assumptions.
To state our main theorem we prepare some notations.

Let $Ksp(\tT^n)$ be the Grothendieck group of
symplectic vector bundles over $\tT^n$,
so we consider the difference $[V_1]- [W_1]$ as
an element in $Ksp(\tT^n)$.
We now introduce a set of invariants to characterize it.

For each subset $S \subset \{1, \cdots, n\}$
we put $\tilde \R^S=\Map(S, \tilde \R)$ 
and induce the Real structure on $\tilde \R^S$ from $\tilde \R$.
Then $\tilde \R^S$ is isomorphic to $\tilde \R^{\caS}$,
where $\caS$ is the the cardinal number of $S$.
We also have the Real torus $\tT^S=\Map(S, \tT^1)$,
which is isomorphic to $\tT^{\caS}$.
Let $\pi_S: \tT^n \to \tT^S$ be the projection
and $i_S: \tilde \R^S \to \tT^S$ an open $\{\pm 1\}$-equivariant
embedding onto an open neighborhood of the origin.
Then, as we shall show later,
the following map is an isomorphism: 
\begin{equation*}
\sum_{S \subset \{1, \cdots , n \}}\pi_S^*(i_S)_!
: \bigoplus_{S}Ksp(\tilde \R^S)
\to Ksp(\tT^n).
\end{equation*}
Thus $[V_1]- [W_1] \in Ksp(\tT^n)$ is determined
by the corresponding components
$a_S \in Ksp(\tilde \R^S)$.
By the Bott periodicity theorem in $Ksp$-theory \cite{Dupont}
(see also Section \ref{real-symplectic})
we see $Ksp(\tilde \R^S) \cong \Z$ if $|S| \equiv 0, 4 \pmod 8$,
$Ksp(\tilde \R^S) \cong \Z/2$ if $|S| \equiv 2$ or $3 \pmod 8$, 
and in the other cases $Ksp(\tilde \R^S)=\{0\}$.

In this paper we investigate the case $a_S=0$ for any nonempty set $S$
with $\caF \equiv 0 \mod 8$. From now we assume it.
Then the quantity $a_S$ for $S$ with  $|S| \equiv 0, 4 \pmod 8$,
which we consider as an integer,
can be detected by the Chern character:
\begin{align}\label{Chern}
\ch(V_1)- \ch(W_1)=2k + \sum_{S \subset \{1, \cdots, n\},\caF=4}
a_S \bigwedge_{i\in S}d\xi^i \in H^*(\tT^n ;\Z),
\end{align}
where $(\xi_1, \cdots, \xi_n)$ is the coordinate of $\tilde \R^n$.

Let $\overline a_S$ be the mod $2$ reduction of $a_S$,
i.e. $\overline a_S =a_S \pmod 2 \in Ksp(\tilde \R^S)\otimes \Z/2 \cong \Z/2$.
Let $X_1, \cdots, X_n$ be indeterminates.
For each $S \subset \{1, \cdots, n \}$
we define an integer $N_S$  by the equation
\begin{align*}
\prod_{\caF=2,3,4, \, \overline a_S \ne 0}(1-2\prod_{i\in S}X_i)
& \equiv \sum_{S \subset \{1, \cdots, n \}}N_S \prod_{i \in S}X_i
\mod {\cal I},
\end{align*}
where $\cal I$ is the ideal generated by
$X_i^2-X_i$ $(1 \le i \le n)$.
Combinatorially $N_S$ is the sum
\begin{align*}
& N_S
= \sum_{m \ge 0}N(S,m)(-2)^m,
\end{align*}
where $N(S,m)$ is the cardinal number of the set
\begin{align*}
\{ \{S_1, \cdots, S_m \}
\mid \,
& \text{$\caFi =2,3$ or $4$, 
$\overline a_{S_i} \ne 0$ ($1 \le i \le m$),
$S= S_1 \cup \cdots \cup S_m$} \}.
\end{align*}
For $S$ with $N_S \ne 0$ we denote by $d_S$ the maximal power
in $2$ which divides $N_S$.
\begin{theorem}\label{main1}
Assume $l > 0$ and $a_S=0$ for any nonempty set $S$ with
$\caF \equiv 0 \mod 8$.
Suppose there is a proper $\pin$-equivariant
fiber-preserving map from $V_0 \oplus V_1$ to $W_0 \oplus W_1$
which induces the identity on the base space $\tT^n$
and whose restriction to $V_0$ is given by
the standard linear inclusion $\tilde \R^x \to \tilde \R^{x+l}$.
If $S \subset \{1, \cdots, n\}$  satisfies $\caS=$ even, $N_S \ne 0$
then we have
\begin{align*}
l\ge 2k+\caS-2d_S+\varepsilon(k+d_S, l), 
\end{align*}
where
\begin{align*}
\varepsilon(\tilde d, l)=
\begin{cases}
3, &  \tilde d \equiv 0 \mod 4, \quad l \ge 4,\\
1, &  \tilde d \equiv 0 \mod 4, \quad l < 4, \\
1, &  \tilde d \equiv 1 \mod 4, \\
2, &  \tilde d \equiv 2 \mod 4, \\
3, &  \tilde d \equiv 3 \mod 4. \\
\end{cases}
\end{align*}
\end{theorem}
We apply Theorem~\ref{main1} to Seiberg-Witten theory \cite{Witten}
by using the argument in \cite{Furuta}.
Then we have the following theorem,
which is our motivation.
The detail of the proof will be given in Section \ref{Proof-of-2}.

Let $X$ be a connected closed oriented spin $4$-manifold
with indefinite intersection form.
Let $\Ind \Bbb D$ be the index bundle
of Dirac operators parameterized by the Jacobian torus
$J_X=H^1(X; \R)/ H^1(X; \Z)$. 
We induce an involution on $J_X$ from
the linear involution on $H^1(X; \R)$
and define a symplectic action on $\Ind \Bbb D$
by the simultaneous action of the involution on $J_X$
and scalar multiplication by $j \in \HH$ on the spinor bundle of $X$.
Thus we consider $\Ind \Bbb D$ as an element in $Ksp(\tT^n)$.
Take a basis $x_1, \cdots, x_n$ of $H^1(X;\Z)$
($n=b_1(X)$) to identify $J_X$ with the Real torus $\tT^n$.
From the above discussion we then obtain the invariant
$a_S \in  KSp(\tilde \R^S)$ for each $S \subset \{1, \cdots,n \}$.

When $\caF=4$ the cohomological formula of the index theorem
\cite{A-S-IV} shows
\begin{align}\label{cup-product}
a_S =\pm\langle \prod_{i \in S}x_i, [X] \rangle.
\end{align}

When $\caF=2$ or $3$ we can describe $a_S$ by using
the the Albanese map $\rho: X \to \hT^n=\Hom(H^1(X; \Z), \R/\Z)$:
The composition of the maps
\begin{align}\label{parity}
\Z \cong
KO^{\caF}(\R^S)
& \overset{(i_S)_!}\to
KO^{\caF}(T^S)
\overset{\pi_S^*}\to KO^{\caF}(\hT^n) \notag\\
& \overset{\rho^*}\to KO^{\caF}(X) 
\overset{i_!}\to KO^{\caF + 4}(\pt)\cong \Z/2
\end{align}
is zero if and only if $\overline a_S=0$,
where $i: X \to \pt$ is the constant map.

\begin{theorem}\label{main2}
Let $X$ be a connected closed oriented spin $4$-manifold
with indefinite intersection form.
Let $n=b_1(X)$.
Then for each $S \subset \{1, \cdots, n \}$
with $\caS$ even and $N_S \ne 0$,
we have
$$b_2^{+}(X) \ge -\dfrac{\sign(X)}8+\caS-2d_S 
+ \varepsilon\left(-\dfrac{\sign(X)}{16}+d_S , b_2^+(X) \right ).$$
\end{theorem}

We take $S =\emptyset$, then $N_{\emptyset}=1$, so we get the corollary below:
\begin{corollary}\label{cor}
Let $X$ be a connected closed oriented spin $4$-manifold
with indefinite intersection form.
Then we have
$$
b_2^{+}(X)
\ge -\dfrac{\sign(X)}8
+ \varepsilon\left(-\dfrac{\sign(X)}{16}, b_2^+(X) \right ).
$$
\end{corollary}

\begin{remark}
Corollary~\ref{cor}
is a small improvement of the inequality obtained in \cite{Furuta}.
Such an improvement began in \cite{FKMa}.
In the case when $-\sign(X)/16 \equiv 2, 3 \mod 4$,
Corollary~\ref{cor} was first shown by N.~Minami \cite{Minami}
and B.~Schmidt \cite{Schmit} independently,
in which they destabilize a $\pin$-equivariant map
to appeal to a result by S.~Stolz in \cite{Stolz},
while our method is more direct.
On the other hand, 
inspired by D.~Ruberman and S.~Strle \cite{Ruberman},
we discussed its improvement when $b_1(X)>0$ in \cite{FKMM}.
Theorem \ref{cor2}, which will be given below,
includes both improvements.
\end{remark}

To get a stronger inequality than Corollary \ref{cor}
we need to compute $d_S$ for a nonempty subset $S$,
which includes calculation of some maps in $KO$-theory.
As discussed in \cite{A-S-V},
this cannot be detected by ordinary cohomology theory.
However we suggest a special case when we can avoid
this calculation.

Let $m T^4$ be the connected sum of $m$-copies of the 
$4$-dimensional torus $T^4=(\R/\Z)^4$.
If we take the standard basis of $H^4(T^4;\Z)$,
it is easy to see $\rank H^4(mT^4;\Z)=4m$
and $N(S_{\max},4m)=1$ for the maximal subset 
$S_{\max}=\{1, \cdots, 4m\}$.
Moreover since $S_{\max}$ cannot be written as
$S_{\max}=S_1 \cup \cdots \cup S_{m'}$ ($2  \le |S_i| \le 4$) 
with $m' < m$,
we see $d_{S_{\max}}=m$
without calculating $a_S$ for $S$ which satisfies $|S|=2$ or $3$.
Hence we get
\begin{theorem}\label{cor2}
Let $X$ be a connected closed oriented spin $4$-manifold
with indefinite intersection form.
If there is an injective homomorphism 
$\iota: H^1(m T^4; \Z) \to H^1(X; \Z)$ which satisfies
$$\langle \iota(x)\iota(y)\iota(z)\iota(w), [X] \rangle
\equiv \langle x y z w, [m T^4] \rangle \mod 2$$
for any $x, y, z, w \in H^1(m T^4; \Z)$,
then we have
$$b_2^{+}(X) \ge -\dfrac{\sign(X)}8 + 2m 
+ \varepsilon\left (-\dfrac{\sign(X)}{16} +m, b_2^+(X) \right ).$$
\end{theorem}

\begin{corollary}\label{cor3}
Let $X$ be a connected closed oriented spin $4$-manifold
with indefinite intersection form.
If $X$ decomposes as
$X =X' \sharp (m T^4)$
by a closed $4$-manifold $X'$
then the inequality
$5m + \varepsilon(-\sign(X)/16+m, b_2^+(X)-3m) 
-\sign(X)/8 \le b_2^{+}(X)$
holds.
\end{corollary}

This paper is organized as follows:
In Section \ref{sec-KO-degree}
we introduce equivariant $KO^*$-degree 
for $\Groupp$-equivariant fiber-preserving maps between 
spin $\Groupp$-vector bundles for a compact Lie group $G$.
In Section \ref{some-properties} we give
some formulae for the Euler class,
which we need in our calculation.
In Section \ref{Products} we discuss how to
calculate the product of elements in $KO_{\Groupp}^*(\pt)$.
In Section \ref{real-symplectic} we recall
Real $K$-theory and symplectic $K$-theory.
Some relations between other $K$-theories are also discussed.
In Section \ref{irreducible}, \ref{calculations}
we carry out calculations in our case.
The proof of Theorem \ref{main1} is given in Section \ref{Proof-of-1},
while we prove  Theorem \ref{main2} in Section \ref{Proof-of-2}.
We also show the equations (\ref{Chern}), (\ref{cup-product}),
(\ref{parity}) in Section \ref{Proof-of-1} and \ref{Proof-of-2}.

\section{Equivariant  $KO^*$-degree}\label{sec-KO-degree}

We first recall the equivariant $KO^*$-theoretic Bott class for spin bundles.

Let $\Groupp$ be a compact Lie group.
For a compact $\Groupp$-space $B$
we let $KO_{\Groupp}(B)$, $K_{\Groupp}(B)$ and $KSp_{\Groupp}(B)$
be the Grothendieck groups of
real, complex and quaternionic $\Groupp$-vector bundles over $B$ respectively.
These definitions extend to a locally compact $\Groupp$-space $B$
in the usual way \cite{Segal}.

Suppose $V$ is a real spin $\Groupp$-vector bundle
over a compact $\Groupp$-space $B$,
i.e. a triple consisting of a real oriented $\Groupp$-vector bundle $V$,
a spin structure $\Spin(V)$ on 
the principal bundle $SO(V)$ of oriented orthonormal frames
(for some $\Groupp$-invariant Riemannian metric on $V$),
and a lift of the $\Groupp$-action on $SO(V)$ to $\Spin(V)$.

If $\rank V =4k$,
the real spinor $\Groupp$-vector bundle $S(V)=S^+(V)\oplus S^-(V)$
is formed from $\Spin(V)$ and
the irreducible $\Z/2$-graded $Cl(\R^{4k})$-representation
$\Delta_{4k}=\Delta^+_{4k}\oplus \Delta^-_{4k}$,
where $Cl(\R^{n})$ is the Clifford algebra (over $\R$) of
the quadratic form $-(x_1^2 +\cdots+ x_n^2)$ on $\R^{n}$ \cite{Atiyah-B-S}.
We follow \cite[Section 6]{Atiyah-B-S} to choose 
the $\Z/2$-grading of $\Delta_{4k}$.
The Clifford multiplication then defines a $\Groupp$-equivariant
bundle map 
$$c: V \to \Hom(S^+(V), S^-(V)).$$
For each vector $v \in V \setminus \{0\}$, $c(v)$ is an isomorphism.

Let $V$ be as above with arbitrary rank.
Let $m$ be a non-negative integer satisfying
$m+\rank V \equiv 0 \bmod 8$,
and $\underbar \R^m$ the product bundle $B \times \R^m$.
Then $V \oplus \underbar\R^m$ is naturally again a spin $\Groupp$-vector bundle
via the inclusion $\Spin(V) \subset \Spin(V \oplus \underbar\R^m)$.
The Bott class $\beta(V \oplus \underbar\R^m)$ is then
an element of $KO_{\Groupp}(V \oplus \underbar\R^m)=KO_{\Groupp}^{-m}(V)$
and defined by the triple
\begin{eqnarray} \label{thom}
[(\pi^* S^+(V \oplus \underbar\R^m),
\pi^* S^-(V \oplus \underbar\R^m),c)].
\end{eqnarray}
where $\pi: V\to B$ is the projection.
By using the Bott periodicity theorem (or the Thom isomorphism theorem)
in equivariant $KO$-theory
\cite{Atiyah},
we have the definition of $KO_{\Groupp}^*$ for positive $*$ and
$KO_{\Groupp}^{-m}(V)$ is identified with $KO_{\Groupp}^{\rank V}(V)$
by this definition. 
With this identification we write
\begin{align*}
\beta(V) \in KO^{\rank V}_{\Groupp}(V)
\end{align*}
for $\beta(V \oplus \underbar\R^m)$.
Then the Bott periodicity theorem 
asserts $KO_{\Groupp}^*(V)$ is a free $KO_{\Groupp}^*(B)$-module
generated by $\beta(V)$.
By the restriction of $\beta(V)$ to $\{0\} \oplus \underbar\R^m$,
we have the Euler class
$$
e(V) \in KO_{\Groupp}^{\rank V}(B).
$$

Let $V$ be again as above with arbitrary rank.
We next use a quaternionic multiplication on $\Delta^{8k+4}$
compatible with the action of $Cl(\R^{8k+4})$,
so we take 
$m'$ to satisfy $m'+\rank V \equiv 4 \bmod 8$.
Then the bundle $S^+(V \oplus \R^{m'})\oplus S^-(V \oplus \R^{m'})$
has a quaternionic $\Groupp$-structure
and the Clifford multiplication $c(v)$ for each vector
$v \in V \setminus \{0\}$ is quaternionic linear.
Thus the $\Groupp$-equivariant $KSp^*$-theoretic Bott class 
$\beta_{\HH}(V)\in KSp_{\Groupp}^{\rank V+4}(V)$
is also defined by the triple (\ref{thom}) 
together with quaternionic scalar multiplication,
if we replace $m$ by $m'$.
Note that multiplication by $\beta_{\HH}(\R^4)$ 
sends $\beta(V)$ to $\beta_{\HH}(V)$.
However we use the explicit construction of $\beta_{\HH}(V)$ later.

Now we define $KO_{\Groupp}^*$-degree for a
proper $\Groupp$-equivariant fiber-preserving map between
spin $\Groupp$-vector bundles.

Suppose $V$ and $W$ are two spin $\Groupp$-vector bundles over 
a compact $\Groupp$-space $B$.
Let $\varphi: V \to W$ be
a proper $\Groupp$-equivariant fiber-preserving map
which induces the identity on the base space $B$.
The $KO_{\Groupp}^*$-degree $\alpha_\varphi$ of $\varphi$ is then
the element of $KO_{\Groupp}^{\rank W -\rank V}(B)$
defined by the relation
\begin{align}\label{KO-degree}
\alpha_\varphi \beta(V) = \varphi^* \beta(W).
\end{align}

Let $\Groupp_0$ be a closed normal subgroup of $\Groupp$
such that $\Groupp_0$ acts trivially on $B$.
Let $V_0$, $W_0$ be spin $\Groupp/\Groupp_0$-vector bundles,
which we also regard spin $\Groupp$-vector bundles.
Let $V_1$, $W_1$ be spin $\Groupp$-vector bundles
whose actions of $\Groupp_0$ are free outside the zeros.

Suppose we have decompositions $V=V_0 \oplus V_1$
and $W=W_0 \oplus W_1$ as spin $\Groupp$-vector bundles.
Then the fixed point sets of $V$ and $W$ for the $\Groupp_0$-actions 
are respectively $V_0$ and $W_0$,
so we have a proper $\Groupp/\Groupp_0$-equivariant fiber-preserving map 
$\varphi_0: V_0 \to W_0$ as a restriction of $\varphi$.
Let $\alpha_{\varphi_0}$ be the $KO^*_{\Groupp/\Groupp_0}$-degree of $\varphi_0$:
\begin{align}\label{KO-degree0}
\alpha_{\varphi_0} \beta(V_0) &= \varphi_0^* \beta(W_0).
\end{align}
We regard $KO_{\Groupp}^*(B)$  as a $KO_{\Groupp/\Groupp_0}^*(B)$-module.
We immediately obtain the following relations
\cite{AtiyahTall}, \cite{Komiya}.
\begin{lemma}\label{keylemma}
\begin{align*}
\alpha_{\varphi_0} e(V_0) &= e(W_0) \in KO_{\Groupp/\Groupp_0}^*(B), \\
\alpha_{\varphi} e(V_1) &= e(W_1)\alpha_{\varphi_0} \in KO_{\Groupp}^*(B).
\end{align*}
\end{lemma}
\proof Restrict \eqref{KO-degree0} to $B$ to
obtain the first relation.
Restrict \eqref{KO-degree} to $V_0$ and use
\eqref{KO-degree0} to obtain the second relation.
\qed

We want to apply this lemma to calculate the $KO_{\Groupp}^*$-degree.
To do so we need to determine the Euler class and 
the product of elements in $KO_{\Groupp}^*(B)$.
In the next two sections we make preparations for it.

\section{Some properties of the Euler class}\label{some-properties}
This section consists of two observations on the Euler class.
First we discuss spin structures on complex or quaternionic
$\Groupp$-vector bundles and the Euler classes for them.
Second we give a formula to calculate the Euler classes
for complex $\Groupp$-vector bundles
coming from the Bott class for a complex representation  of $\Groupp$.
From now on we deal with operations between
various equivariant $K$-theories.
Thus we begin by introducing notation.
The following operations are elementary:

\begin{enumerate}
\item  Let $V$ be a real ${\Groupp}$-vector bundle over $B$.
Denote by $cV=\C \otimes_{\R} V$ the complexification of $V$.
This operation induces a homomorphism $c: KO_{\Groupp}(B) \to K_{\Groupp}(B)$.
\item  Let $V$ be a complex ${\Groupp}$-vector bundle over $B$.
Denote by $qV=\HH \otimes_{\C} V$ the quaternization of $V$.
This operation induces a homomorphism $q: K_{\Groupp}(B) \to KSp_{\Groupp}(B)$.
\item  Let $V$ be a complex ${\Groupp}$-vector bundle over $B$.
Denote by $rV$ the restriction of the scalars from $\C$ to $\R$.
We denote by $r: K_{\Groupp}(B) \to KO_{\Groupp}(B)$
the induced homomorphism.
\item  Let $V$ be a quaternionic ${\Groupp}$-vector bundle over $B$.
Denote by $c'V$ the restriction of the scalars from $\HH$ to $\C$.
We denote by $c': KSp_{\Groupp}(B) \to K_{\Groupp}(B)$
the induced homomorphism.
\item  Let $V$ be a complex ${\Groupp}$-vector bundle over $B$.
Denote by $tV=\overline V$ the complex conjugate of $V$.
This operation induces a homomorphism $t: K_{\Groupp}(B) \to K_{\Groupp}(B)$.
\item  Let $V$ be a complex ${\Groupp}$-vector bundle over $B$.
Since the action of $\Groupp$ commutes with the complex multiplication,
$V$ is naturally  regarded as a complex ${\Groupp}\times U(1)$-vector bundle.
We denote this complex ${\Groupp}\times U(1)$-vector bundle by $bV$.
Let $b: K_{\Groupp}(B) \to K_{\Groupp \times U(1)}(B)$ be 
the induced homomorphism.
\item Let $V$ be a quaternionic $\Groupp$-vector bundle $V$ over $B$.
Since the action of $\Groupp$ commutes with the quaternionic multiplication,
$rc'V$ is naturally regarded as a real $\Groupp \times Sp(1)$-vector bundle.
We denote this real $\Groupp \times Sp(1)$-vector bundle by $\rc V$.
Let $i: U(1) \to Sp(1)$  be the inclusion.
Then it is obvious $i^*\rc V=rbc' V$.
We let $s: KSp_{\Groupp}(B) \to KO_{\Groupp \times Sp(1)}(B)$
be the induced homomorphism.
\end{enumerate}

Suppose $V$ is a complex $\Groupp$-vector bundle
over a compact $\Groupp$-space $B$.
Let $\Lambda^*(V)=\Lambda^{\even}(V) \oplus \Lambda^{\odd}(V)$
be the exterior power of $V$.
Then in this paper the $K_{\Groupp}$-theoretic Bott class 
$\beta_{\C}(V) \in K_{\Groupp}(V)$
is the triple 
$$[(\pi^*\Lambda^{\even}(V), \pi^*\Lambda^{\odd}(V), \alpha)],$$
where the bundle map
$$\alpha: V \to \Hom(\Lambda^{\even}(V), \Lambda^{\odd}(V))$$
is given by $\alpha(v)= v \wedge \cdot - v \lrcorner \cdot $
\cite{Atiyah-B-S}.
The Euler class 
$e_{\C}(V) \in K_{\Groupp}(B)$
is the restriction to the zero.

Let $\det V$ denote the determinant bundle
of a complex $\Groupp$-vector bundle $V$.
If  $\det V$ has a square root $\det(V)^{1/2}$
as a complex $\Groupp$-line bundle
then the complex  $\Groupp$-structure of $V$ 
canonically induces a spin $\Groupp$-structure on 
the real $\Groupp$-vector bundle $rV$ (see \cite{Atiyah-B-S}).
Then the Bott class $\beta_{\C}(V)$
and the complexification $c(\beta(rV))$ of the Bott class $\beta(rV)$
are related by
\begin{align}\label{Bott-c}
c(\beta(rV))=\det(V)^{-1/2} \beta_{\C}(V)\in K_{\Groupp}(V),
\end{align}
so we have 
\begin{align}\label{Euler-c}
c(e(rV))=\det(V)^{-1/2} e_{\C}(V)\in K_{\Groupp}(B).
\end{align}

The above argument is obviously applicable to quaternionic
$\Groupp$-vector bundles.
However it is insufficient in our application,
so we will give its refinement.

Let $V$ be a quaternionic $\Groupp$-vector bundle
over a compact $\Groupp$-space $B$ with $\rank_{\HH} V =n$.
We first show that
there is a natural spin $\Groupp\times Sp(1)$-structure on 
the real $\Groupp \times Sp(1)$-vector bundle $\rc V$.
Let $Sp(V)$ be the principal $Sp(n)$-bundle of orthonormal frames
for some $\Groupp$-invariant quaternionic metric on $V$.
The principal $SO(4n)$-bundle $SO(\rc V)$
is then $\Groupp \times Sp(1)$-isomorphic to 
$(Sp(V) \times Sp(1) \times SO(4n))/(Sp(n)\times Sp(1))$,
where the action of $Sp(n)\times Sp(1)$ on $SO(4n)$
is factoring through the homomorphism
$$Sp(n)\times Sp(1) \to (Sp(n)\times Sp(1))/\{\pm 1\}
\to SO(4n).$$
The lift to $\Spin(4n)$ 
provides a spin $\Groupp\times Sp(1)$-structure on
$\rc V$.

On the other hand 
there is also a natural spin $\Groupp\times U(1)$-structure on 
the real $\Groupp\times U(1)$-vector bundle on $rbc'V$,
since the principal $SO(4n)$-bundle
$SO(rbc'V)$ is $\Groupp \times U(1)$-isomorphic to
$(Sp(V)\times U(1) \times SO(4n))/(Sp(n)\times U(1))$
and the action of $Sp(n)\times U(1)$ on $SO(4n)$
is factoring through the homomorphism
$$Sp(n)\times U(1) \to SU(2n)\times U(1) \to 
(SU(2n)\times U(1))/\{\pm 1\}\to SO(4n).$$

Let $i: U(1) \to Sp(1)$ be the inclusion.
Then from the above constructions
the spin $\Groupp \times U(1)$-structure
of $i^* \rc V$ is obviously isomorphic to that of $rbc'V$.
Moreover under the canonical trivialization
$\det c'V \cong B \times \C_{2n}(t)$
the spin $\Groupp \times U(1)$-structure
of $c'V$ is isomorphic to the one obtained from the square root
$(\det c'V)^{1/2}=B \times \C_n(t)$, 
where $\C_n(t)$ is the one-dimensional
complex representation of $U(1)$ with weight $n$.
Hence we use (\ref{Euler-c}) to get
\begin{proposition}\label{Euler-class-0}
Let $V$ be a quaternionic $\Groupp$-vector bundle over $B$
with $\rank_{\HH} V=n$.
Then we have 
$$
c(e(i^*\rc V))
=\C_{-n}(t) e_{\C}(bc'V) \in K_{\Groupp \times U(1)}(B),
$$
where $i: U(1) \to Sp(1)$ is the inclusion.
\end{proposition}

Now we let $V$ be a complex representation of $\Groupp$.
Then the map $\alpha$ gives an isomorphism 
between the product $\Groupp$-vector bundles 
$V \setminus \{0\} \times \Lambda^{\even}(V)$
and $V \setminus \{0\}\times \Lambda^{\odd}(V)$.
If we regard $\alpha$ as a trivialization of
$V \times \Lambda^{\even}(V)$ over $V\setminus \{0\}$
then we obtain a complex $\Groupp$-vector bundle $\hat V^+$ 
as an extension of $V \times \Lambda^{\even}(V)$ to
the one-point compactification $V \cup \{\infty\}$ of $V$, namely the
$\Groupp$-sphere.
If we exchange  $\Lambda^{\even}(V)$ for
$\Lambda^{\odd}(V)$ and replace $\alpha$ by $\alpha^{-1}$,
we obtain another complex $\Groupp$-vector bundle $\hat V^-$
over $V \cup \{\infty\}$.
Under the canonical decomposition
$K_{\Groupp}(V \cup \{\infty\})
\cong K_{\Groupp}(\{\infty\})\oplus K_{\Groupp}(V)$
these $\Groupp$-vector bundles satisfy
\begin{align*}
[\hat V^+] & = [\Lambda^{\odd}(V)]+\beta_{\C}(V),
\quad [\hat V^-]  = [\Lambda^{\even}(V)]-\beta_{\C}(V).
\end{align*}

\begin{definition}
Let $x_0$, $x_1, \cdots, x_n$ be indeterminates.
We define a polynomial $\mu_n(x_0, x_1, \cdots, x_n)$ 
by the identity{\rm :}
\begin{align*}
& \mu_n(x_0, x_1, \cdots, x_n) \prod_{1 \le i \le n}(1-x_i) \\
& =\prod_{S \subset \{1, \cdots, n \}, \caS: \even}
(1-x_0 \prod_{i \in S}x_i)
-\prod_{S \subset \{1, \cdots, n \}, \caS: \odd}
(1-x_0 \prod_{i \in S}x_i).
\end{align*}
For instance $\mu_1=\mu_2=-x_0$ and $\mu_3=-x_0+x_0^3x_1x_2x_3$.
\end{definition}

\begin{proposition}\label{Euler-class-1}
Let $L_1, \cdots, L_n$ be complex $1$-dimensional representations of $\Groupp$
and $V =L_1 \oplus \cdots \oplus L_n$.
Let $L_0$ be another complex $1$-dimensional representation of $\Groupp$. 
Then we have
\begin{align*}
e_{\C}(L_0 \otimes \hat V^+) &=e_{\C}(L_0 \otimes \Lambda^{\odd}(V))
+\mu_n(L_0, L_1, \cdots, L_n)\beta_{\C}(V),\\
e_{\C}(L_0 \otimes \hat V^-)& = e_{\C}(L_0 \otimes \Lambda^{\even}(V))
-\mu_n(L_0, L_1, \cdots, L_n)\beta_{\C}(V).
\end{align*}
\end{proposition}
\proof
We introduce an action of $U(1)^{n+1}$ on $L_i$ 
via the $i$-th projection from $U(1)^{n+1}$ to $U(1)$ 
and multiplication ($0 \le i \le n$).
Since the action of $\Groupp$ on $L_i$ factors through a homomorphism
$\Groupp \to U(1) \to U(1)^{n+1}$, 
it is sufficient to show our formula in $K_{U(1)^{n+1}}(V \cup \{\infty\})$.
By the Bott periodicity theorem
$K_{U(1)^{n+1}}(V \cup \{\infty\})$ is freely generated by 
$1$ and $\beta_{\C}(V)$ as a $R(U(1)^{n+1})$-algebra,
so we can write $e_{\C}(L_0 \otimes \hat V^+)$ as
$$e_{\C}(L_0 \otimes \hat V^+)= \alpha_n + \beta_n \beta_{\C}(V)$$
for some $\alpha_n, \beta_n \in R(U(1)^{n+1})$.
Restricting this equation to $\infty$ and $0$, we get
\begin{align*}
\alpha_n & =e_{\C}(L_0 \otimes \Lambda^{\odd}(V)),\\
\beta_ne_{\C}(V) &=
e_{\C}(L_0 \otimes \Lambda^{\even}(V))-e_{\C}(L_0 \otimes\Lambda^{\odd}(V)).
\end{align*}
Since the representation ring $R(U(1)^{n+1})$ has no zero-divisor,
we have $\beta_n=\mu_n(L_0,L_1, \cdots,L_n)$.
We can compute the Euler class
$e_{\C}(L_0 \otimes \hat V^-)$ in the same way.
\qed

Here we suggest a method to calculate the polynomial 
$\mu_n(x, 1, \cdots, 1)$,
which will be appeared in our application.
Although we shall not need its explicit form in it,
it may be helpful in understanding this polynomial.
\begin{definition}
Let  $\nu_n(x)$ be the polynomial defined by 
$$\nu_n(x)=-(1-x)^{2^{n-1}} s_{n-1}(x), \quad
s_0(x)=x/(1-x), \quad s_n(x)=x s_{n-1}'(x).$$
\end{definition}

\begin{proposition}
The identity $\mu_n(x,1, \cdots, 1)=\nu_n(x)$ holds for any $n$.
\end{proposition}

We give a geometric proof although the identity is purely algebraic.
Then the above proposition is clearly equivalent to
\begin{proposition}\label{Euler-class-2}
Suppose $V$ is the trivial $n$-dimensional representation of $\Groupp$.
Let $L$ be a complex $1$-dimensional representation of $\Groupp$. 
Then we have
\begin{align*}
e_{\C}(L \otimes \hat V^+)& =e_{\C}(L \otimes \Lambda^{\odd}(V))
+\nu_n(L) \beta_{\C}(V),\\
e_{\C}(L \otimes \hat V^-)& =e_{\C}(L \otimes \Lambda^{\even}(V))
-\nu_n(L) \beta_{\C}(V).
\end{align*}
\end{proposition}
\proof
We use some operations in  $K$-theory (cf. \cite{Atiyah-K}).
Let $\lambda_t[L \otimes \hat V^+]
=\sum_{i=0}^n t^i[\Lambda^i(L \otimes \hat V^+)]
\in 1+K_{\Groupp}(V \cup \{\infty\})[[t]]$.
Then we have 
$$
e_{\C}(L \otimes \hat V^+)= \lambda_{-1}[L \otimes \hat V^+]
=(\lambda_{t}[L\otimes \Lambda^{\odd}(V)]\lambda_{t}[L\otimes \beta_{\C}(V)])
\vert_{t=-1}.
$$
Let $\psi^k$ be the Adams operation.
We write $W= L \otimes \beta_{\C}(V)$.
Then $\psi^k(W)=k^n L^k\beta_{\C}(V)$
and the formal power series $\psi_t(W)=
\sum_{k=1}^{\infty}t^k \psi^k(W)$
is computed as
$\psi_t(W)
=s_n(tL) \beta_{\C}(V)
$
by using $s_n(x)=1^n x +2^n x^2 + \cdots + k^nx^k + \cdots$.
By the Newton formula
(cf. Chapter 13 of \cite{Husemoller}),
the differential equation
$$\psi_{-t}(W)
=-t \lambda_{t}'(W)/\lambda_t(W),
\quad\lambda_0(W)=1
$$
determines $\lambda_t(W)$ uniquely,
so we obtain
$\lambda_t(W)=1-s_{n-1}(-tL)\beta_{\C}(V)$.

By Proposition \ref{Euler-class-1}
the coefficient of $\beta_{\C}(V)$  must be $-\nu_n(L)$
for $e_{\C}(L_0 \otimes \hat V^-)$.
\qed

\section{Calculations of the products}\label{Products}

Let $\Groupp$ be a compact Lie group.
For a real irreducible representation space $V$ of $\Groupp$,
the ring $\End^{\Groupp}(V)$ of the $\Groupp$-invariant endomorphisms
is a field. We have the following three cases.
\begin{itemize}
\item
Real case: $\End^{\Groupp}(V) \cong \R$.
We denote by $\Irr_{\R}$  the set consisting of
these representations.
\item
Complex case: $\End^{\Groupp}(V) \cong \C$.
We denote by $\Irr_{\C}$  the set consisting of
these representations.
\item
Quaternionic case:
$\End^{\Groupp}(V) \cong \HH$.
We denote by $\Irr_{\HH}$  the set consisting of
these representations.
\end{itemize}
Then any real representation $W$ of $\Gamma$
can be naturally decomposed as
$$ \bigoplus _{\text{$ V \in \Irr_{\F},\F=\R, \C, \HH$}}
V \otimes_{\End^{\Groupp}(V)} \Hom^{\Groupp}(V, W) \overset{\cong}
\longrightarrow U.$$
Thus this decomposition still holds for a $\Groupp$-vector bundle $W$
if the action of $G$ on the base space $B$ is trivial.
Hence, as shown by Segal \cite{Segal},
if we fix an isomorphism $\End^{\Groupp}(V) \cong \F$
for each $V \in \Irr_{\F}$ and $\F=\R, \C, \HH$,
we have a canonical decomposition:
\begin{align}\label{KOG}
KO^{q}_{\Groupp}(\pt) \cong  \Z\Irr_{\R} \otimes KO^{q}(\pt)
\oplus \Z\Irr_{\C}  \otimes K^{q}(\pt)
\oplus  \Z\Irr_{\HH} \otimes KSp^{q}(\pt),
\end{align}
which extends to the case $q > 0$ by the Bott periodicity theorem.

We should remark that
two different choices of isomorphisms $\End^{\Groupp}(V) \cong \F$ give
different isomorphisms in (\ref{KOG}) only  when $\F=\C$.
From now on we fix an isomorphism 
$\End^{\Groupp}(V) \cong \C$ for each $V \in \Irr_{\C}$.

For our later use we fix generators of $KO^{q}(\pt)$, $K^{q}(\pt)$
and $KSp^{q}(\pt)$ in the following way:
Put
\begin{align*}
\beta_{\R,8k}& =\beta(\R^{8k}) \in KO^{8k}(\pt) \cong \Z,\\
\beta_{\C,2k}& =\beta_{\C}(\C^{2k}) \in K^{2k}(\pt) \cong \Z,\\
\beta_{\HH,8k+4}& =\beta_{\HH}(\R^{8k+4})\in KSp^{8k+4}(\pt) \cong \Z.
\end{align*}

Then $\beta_{\R,8k+4} =rc'\beta_{\HH,8k+4}$
and $\beta_{\HH,8k} =qc\beta_{\R, 8k}$ are also generators of
$KO^{8k+4}(\pt) \cong \Z$ and $KSp^{8k}(\pt) \cong \Z$ respectively.

Let $\beta_{\R,8k+i}\in KO^{8k+i}(\pt)\cong  \Z/2$ ($i=6,7$),
$\beta_{\HH,8k+i}\in KSp^{8k+i}(\pt)\cong\Z/2$ ($i=2,3$)
be the generators.

\begin{remark}\label{relation}
In our convention we have the following relations
(cf. \cite{Atiyah-B-S}).
\begin{alignat*}{2}
c\beta_{\R,8k}& =\beta_{\C,8k},
\quad&  c\beta_{\R,8k+4}& =2\beta_{\C,8k+4},\\
r\beta_{\C,8k}& =2\beta_{\R,8k},
\quad & r\beta_{\C,8k+4}& =\beta_{\R,8k+4},\\
c'\beta_{\HH,8k+4}& =\beta_{\C,8k+4}, 
\quad & c'\beta_{\HH,8k}& =2\beta_{\C,8k},\\
t\beta_{\C,2k}& =(-1)^k\beta_{\C,2k},
& & \\
\beta_{\R,8k+6} & =r\beta_{\C,8k+6} ,
\quad & \beta_{\HH,8k+2} & =q\beta_{\C, 8k+2},\\
\beta_{\HH, 8(k+1)+2}& =\beta_{\R,8k+6}\beta_{\HH, 4},
\quad & \beta_{\HH, 8(k+1)+3}& =\beta_{\R,8k+7}\beta_{\HH, 4},\\
\beta_{\R, 8(k+1)+6}& =\beta_{\R,8k+7}\beta_{\R, 7},
\quad & \beta_{\HH, 8(k+1)+2}& =\beta_{\R,8k+7}\beta_{\HH, 3}.
\end{alignat*}
\end{remark}
Let $V \in \Irr_{\F}$.
If $K\F^{q}(\pt)$ is nontrivial, we define
\begin{align}\label{nontrivial}
[V]_{q}=[V] \, \beta_{\F,q} 
\in \Z\Irr_{\F} \otimes K\F^q(\pt) \subset KO_{\Groupp}^{q}(\pt),
\end{align}
where $K\R=KO$, $K\C=K$ and $K\HH=KSp$.
Then by the decomposition (\ref{KOG})
any element in $KO_{\Groupp}^q(\pt)$
can be uniquely written as a linear combination of $[V]_q$
for $V \in \Irr_{\F}$
with $K{\bf F}^{q}(\pt)$ nontrivial, where $\F=\R, \C$ or $\HH$.
We should note that
\begin{align*}
[V]_{2} &=[V] \,  q\beta_{{\bf \C},2}
=[c'V] \, \beta_{{\bf C},2}
=[c'V]_{2}, \quad \text{if $\F=\HH$},\\
[V]_{6} & = [V] \, r\beta_{{\bf C},6}
= [cV] \, \beta_{{\bf C},6}
=[cV]_{6}, \quad \text{if $\F=\R$}.
\end{align*}

We want to calculate the product
\begin{align*}
KO_{\Groupp}^{p}(\pt) \times KO^{q}_{\Groupp}(\pt) 
& \to KO^{p+q}_{\Groupp}(\pt).
\end{align*}
under the decomposition (\ref{KOG}).
To do so we use the canonical isomorphisms \cite{Adams}:
\begin{align}\label{three}
i+ r+rc'i & : \Z \Irr_{\R} \oplus \Z \Irr_{\C} \oplus \Z \Irr_{\HH}
\to RO({\Groupp}),\notag \\
ci + i + ti + c'i & : \Z \Irr_{\R} \oplus \Z \Irr_{\C} \oplus
\Z \Irr_{\C}\oplus \Z \Irr_{\HH} \to R({\Groupp}),\\
qci +qi + i  & : \Z \Irr_{\R} \oplus \Z \Irr_{\C} \oplus \Z \Irr_{\HH}
\to RSp({\Groupp}), \notag
\end{align}
where $i: \Z \Irr_{\F} \to R\F({\Groupp})$ is
the inclusion and $R\R=RO$, $R\C=R$, $R\HH=RSp$.
Note that we used the fixed isomorphism $\End^{\Groupp}(V)\cong \C$ 
for $V \in \Irr_{\C}$.

\begin{lemma}\label{complexi}
The complexification map $c: KO^{2k}_{\Groupp}(\pt) \to
K^{2k}_{\Groupp}(\pt) \cong R(\Groupp)$ is given by
\begin{align*}
c[V]_{2k}& = [V] +(-1)^k[tV], \quad V \in \Irr_{\C},\\
c[V]_{4}&=
\begin{cases}
2[cV], & V \in \Irr_{\R}, \\
[c'V],  & V \in \Irr_{\HH}.
\end{cases}
\end{align*}
Moreover the kernel of  $c$ coincides with the torsion subgroup.
\end{lemma}
\proof
Let $V \in \Irr_{\F}$.
When $\F =\R$ or $\C$, 
this easily follows from Remark \ref{relation}.
When $\F=\HH$, we use
the fact that
$\C \otimes_{\R} \Delta_4 \cong \Lambda^*(\C^2) \otimes _{\C}\HH$
as $\Z/2$-graded $\C \otimes_{\R}Cl(\R^4)\otimes_{\R} \HH$-modules.
(Compare the complex $\Spin(4) \times Sp(1)$-modules coming from
the inclusions $\Spin(4) \subset Cl(\R^4)$ and $Sp(1) \subset \HH$.)
Then we have $c([V] \, \beta_{\HH, 4})=[c'V] \, \beta_{\C, 4} $.
The second isomorphism in (\ref{three}) guarantees the kernel of $c$
is in the torsion subgroup.
\qed

This lemma asserts
when $p+q \equiv 0 \pmod 4$
the product can be computed from the complex representation ring of
$\Groupp$.
The other cases can be computed as follows:

We first note that the definition (\ref{nontrivial})
obviously extends to any element $[V] \in R\F(\Groupp)$
by using the tensor product over $\F$,
thus we get an element $[V]_q \in KO_{\Groupp}^q(\pt)$,
if $K\F^{q}(\pt)$ is nontrivial.
To decompose $[V]_q$ according to (\ref{KOG})
we again use the isomorphisms (\ref{three}).
\begin{lemma}\label{decomposed}
Let $[V] \in R\F(\Groupp)$ be decomposed as
$$
[V]=
\begin{cases}
V_{\R} + rV_{\C} + rc'V_{\HH}, & \F=\R,\\
cV_{\R} + V_{\C} + tV'_{\C}+ c'V_{\HH}, & \F=\C,\\
qcV_{\R} + qV_{\C} + V_{\HH}, & \F=\HH,
\end{cases}
$$
where $V_{\R} \in \Z \Irr_{\R}$, $V_{\C}, V'_{\C} \in \Z \Irr_{\C}$,
$V_{\HH} \in \Z \Irr_{\HH}$.
Then we have
\begin{align*}
[V]_{2} & =
\begin{cases}
[V_{\C}]_{2} -[V'_{\C}]_{2}+[V_{\HH}]_{2},
& \F=\C, \\
[V_{\HH}]_{2}, 
& \F=\HH,
\end{cases}\\
[V]_{3} & =[V_{\HH}]_{3}, \quad \F=\HH, \\ 
[V]_{6} & =
\begin{cases}
[V_{\R}]_6,
& \F=\R, \\
[V_{\R}]_6+[V_{\C}]_6 -[V'_{\C}]_6, 
& \F=\C, \\
\end{cases}\\
[V]_{7} & =[V_{\RR}]_{7}, \quad \F=\R. 
\end{align*}
\end{lemma}
\proof
From Remark \ref{relation} we see
\begin{align*}
[tV'_{\C}] \, \beta_{\C, 2k} 
& =t([tV'_{\C}] \, \beta_{\C, 2k} )
=(-1)^k[V'_{\C}]_{2k} \in KO_{\Groupp}^{2k}(\pt), \\
[c'V_{\HH}] \, \beta_{\C, 2}
& = [V_{\HH}] \, q\beta_{\C, 2}
=[V_{\HH}]_2\in KO_{\Groupp}^{2}(\pt),\\
[cV_{\R}] \, \beta_{\C, 6}
& = [V_{\R}] \, r\beta_{\C, 6}  =[V_{\R}]_6\in KO_{\Groupp}^{6}(\pt).
\end{align*}
Moreover since $c': Ksp^{2}(\pt) \to K^2(\pt)$,
$c: KO^{6}(\pt) \to K^6(\pt)$ are obviously the zero maps,
we have
\begin{align*}
[qV_{\C}] \, \beta_{\HH, 2} 
& = [V_{\C}] \,  c'\beta_{\HH,2} =0,\\
[rV_{\C}] \, \beta_{\R, 6}
& = [V_{\C}] \, c\beta_{\R,6} =0.
\end{align*}
Then the identities for $[V]_2$ and $[V]_6$ immediately follow. 
The other cases are trivial.
\qed

\begin{lemma}\label{product}
Let $V_0\in \Irr_{\F_0}$ and $V_1 \in \Irr_{\F_1}$.
Then the products $[V_0]_0[V_1]_{2}$, $[V_0]_0[V_1]_{6}$,
$[V_0]_4[V_1]_{2}$ and $[V_0]_4[V_1]_{6}$
are given by the tables below,
where the symbol $\otimes$ in the tables means the tensor product over $\C$.
\begin{table}[h]\label{table1}
\caption{$[V_0]_0[V_1]_{2}$}
\begin{center}
\begin{tabular}{c|cc}
&  $\F_1=\C$ & $\F_1=\HH$ \\
\hline
$\F_0=\R$
& $[cV_0 \otimes V_1]_{2}$
& $[cV_0 \otimes  c'V_1]_{2}$ \\
$\F_0=\C$
& $[V_0 \otimes  V_1]_{2}$ 
& $[V_0 \otimes  c'V_1]_{2}$\\ 
$\F_0=\HH$
& $[c'V_0 \otimes  V_1]_{2}$
& $[c'V_0 \otimes c'V_1]_{2}$
\end{tabular}  
\end{center} 
\end{table}
\begin{table}[h]\label{table2}
\caption{$[V_0]_0[V_1]_{6}$}
\begin{center}
\begin{tabular}{c|cc}
&  $\F_1=\R$ & $\F_1=\C$ \\
\hline
$\F_0=\R$
& $[cV_0 \otimes cV_1]_{6}$ 
& $[cV_0 \otimes V_1]_{6}$\\
$\F_0=\C$
& $[V_0 \otimes  cV_1]_{6}$
& $[V_0 \otimes  V_1]_{6}$ \\
$\F_0=\HH$
& $[c'V_0 \otimes cV_1]_{6}$ 
& $[c'V_0 \otimes  V_1]_{6}$\\
\end{tabular}  
\end{center} 
\end{table}
\begin{table}[h]\label{table3}
\caption{$[V_0]_4[V_1]_{2}$}
\begin{center}
\begin{tabular}{c|ccc}
&   $\F_1=\C$ & $\F_1=\HH$ \\
\hline
$\F_0=\R$
& $2[cV_0 \otimes V_1]_{6}$ 
& $0$ \\
$\F_0=\C$
& $[V_0 \otimes  V_1]_{6}  +[tV_0 \otimes  V_1]_{6}$
& $0$ \\
$\F_0=\HH$
& $[c'V_0 \otimes V_1]_{6}$
& $[c'V_0 \otimes c'V_1]_{6}$
\end{tabular}  
\end{center} 
\end{table}
\begin{table}[h]\label{table4}
\caption{$[V_0]_4[V_1]_{6}$}
\begin{center}
\begin{tabular}{c|ccc}
&  $\F_1=\R$ & $\F_1=\C$  \\
\hline
$\F_0=\R$
& $0$
& $2[cV_0 \otimes V_1]_{2}$ \\
$\F_0=\C$
& $0$
& $[V_0 \otimes  V_1]_{2} +[tV_0 \otimes  V_1]_{2}$\\
$\F_0=\HH$
& $[c'V_0 \otimes  cV_1]_{2}$
& $[c'V_0 \otimes V_1]_{2}$\\
\end{tabular}  
\end{center} 
\end{table}

\end{lemma}
\proof
Let $\F_0=\F_1=\R$.
Then 
$$
[V_0]_0[V_1]_{6}
=[V_0 \otimes_{\R}V_1] \, r\beta_{\C,6}
=[cV_0 \otimes_{\C}cV_1] \, \beta_{\C,6}
=[cV_0 \otimes_{\C}cV_1]_6.
$$ 
The other cases in the tables of $[V_0]_0[V_1]_{2}$ 
and $[V_0]_0[V_1]_{6}$ 
also follow from the relations $\beta_{\R, 6}=r\beta_{\C,6}$ and 
$\beta_{\HH,2}=q\beta_{\C,2}$ in Remark \ref{relation}. 

We next consider the tables of $[V_0]_4[V_1]_{2}$ and $[V_0]_4[V_1]_{6}$.
Now fix the canonical isomorphism
$Cl(\R^4) \otimes_{\R} Cl(\R^{4k+2}) \cong Cl(\R^{4k+6})$
as $\Z/2$-graded $\R$-algebras \cite{Atiyah-B-S}.

Let $\F_0=\R$.
Counting dimension shows
\begin{align*}
\Delta_4\otimes_{\R}\Lambda^*(\C^{2k+1})
\cong 
\Lambda^*(\C^{2k+3})\oplus\Lambda^*(\C^{2k+3})
\end{align*}
as $\Z/2$-graded $Cl(\R^{4k+6})\otimes_{\R}\C$-modules.
Thus, for instance if $\F_1=\R$, we obtain
$$
[V_0]_4[V_1]_{6}
=[cV_0 \otimes_{\C}V_1] \, 2r\beta_{\C,10}
=0.
$$ 

Let $\F_0=\HH$.
If we consider $\Delta_4$
as a $\Z/2$-graded quaternionic $Cl(\R^4)$-module, we see
\begin{align*}
\Delta_4\otimes_{\R}\Lambda^*(\C^{2k+1})
\cong 
\HH \otimes_{\C}\Lambda^*(\C^{2k+3})
\end{align*}
as $\Z/2$-graded $\HH \otimes_{\R} Cl(\R^{4k+6})\otimes_{\R}\C$-modules.
(Compare the complex $Sp(1) \times \Spin(4k+6)$-modules
coming from the inclusions
$Sp(1) \subset \HH$ and $\Spin(4k+6) \subset Cl(\R^{4k+6})$.)
Thus, for instance if $\F_1=\R$, we obtain
$$
[V_0]_4[V_1]_{6}
=[V_0 \otimes_{\R}V_1] \, q\beta_{\C,10}
=[c'V_0 \otimes_{\C}cV_1] \, \beta_{\C,10}
=[c'V_0 \otimes_{\C}cV_1]_{10},
$$ 
where we induce an $\HH$-action on $V_0 \otimes_{\R}V_1$ 
from that on $V_0$.

Let $\F_0=\C$.
In this case we have
\begin{align*}
& \Lambda^*((\C')^2)
\otimes_{\R}\Lambda^*(\C^{2k+1})
\cong 
\C'\otimes_{\R}\Lambda^*(\R^2)
\otimes_{\R}\Lambda^*(\C^{2k+1})
\\
& \cong \C'\otimes_{\R}\Lambda^*(\C^{2k+3})
\cong 
\Lambda^*(\C^{2k+3}) \oplus 
t\Lambda^*(\C^{2k+3})
\end{align*}
as $\Z/2$-graded $\C' \otimes_{\R} Cl(\R^{4k+6})\otimes_{\R}\C$-modules,
where $\C'$ is a copy of $\C$,
and where $\C$ acts on $t\Lambda^*(\C^{2k+3})$ via scalar 
multiplication on $\C$,
while $\C'$ acts on $t\Lambda^*(\C^{2k+3})$ via complex conjugation.
Thus if we identify $\C'$ with $\C$ by the complex conjugation
then the map
$$\alpha: (\C')^2 \oplus \C^{2k+1}
\to \Hom(\Lambda^{\even}(\C^{2k+3}), \Lambda^{\odd}(\C^{2k+3}))$$
defined in Section \ref{sec-KO-degree} is complex linear.
However the complex conjugation on $(\C')^2$
is homotopic to the identity in $SO((\C')^2)$.
(This fact can be also used to show $t\beta_{\C,-4}=\beta_{\C,-4}$.)
Hence, for instance if $\F_1=\C$, we have
\begin{align*}
[V_0]_{4}[V_1]_{4k+2}
& = [V_0 \otimes_{\C}V_1] \, \beta_{\C,4k+6}
+[tV_0 \otimes_{\C}V_1] \, \beta_{\C,4k+6}\\
&= [V_0 \otimes_{\C}V_1]_{4k+6}
+[tV_0 \otimes_{\C}V_1]_{4k+6}.
\end{align*}
The other cases also follow
from the relations $\beta_{\R, 6}=r\beta_{\C,6}$ and 
$\beta_{\HH,2}=q\beta_{\C,2}$.
\qed

The only remaining nontrivial case to be calculated is
\begin{align*}
KO_{\Groupp}^{4k+3}(\pt) \times KO^{q}_{\Groupp}(\pt) 
& \to KO^{4k+3+q}_{\Groupp}(\pt).
\end{align*}
\begin{lemma}
Let $V_0\in \Irr_{\F_0}$ and $V_1 \in \Irr_{\F_1}$.
Then the products $[V_0]_{3}[V_1]_{q}$ and 
$[V_0]_{7}[V_1]_{q}$
are given by the tables below,
where if $\F_0=\F_1=\HH$ then $V_0 \otimes_{\HH} V_1$
is regarded as a real representation of $\Groupp$,
and if $(\F_0, \F_1)=(\HH, \R)$ or $(\F_0, \F_1)=(\R, \HH)$
then $V_0 \otimes_{\R} V_1$
is regarded as a quaternionic  representation of $\Groupp$.
\begin{table}[h]\label{table5}
\caption{$[V_0]_{3}[V_1]_{q}$}
\begin{center}
\begin{tabular}{c|ccc}
$\F_0=\HH$ &  $\F_1=\R$ & $\F_1=\C$ & $\F_1=\HH$ \\
\hline
$q=0$
& $[V_0 \otimes_{\R} V_1]_{3}$ 
& $[V_0 \otimes_{\R} rV_1]_{3}$
& $[V_0 \otimes_{\R}  rc'V_1]_{3}$ \\
$q=3$
& $0$
& $0$ 
& $[V_0 \otimes_{\HH} V_1]_{6}$\\
$q=4$
& $0$ 
& $0$
& $[V_0 \otimes_{\HH} V_1]_{7}$\\
$q=7$
& $[V_0 \otimes_{\R} V_1]_{2}$ 
& $0$
& $0$
\end{tabular}  
\end{center} 
\end{table}
\begin{table}[h]\label{table6}
\caption{$[V_0]_{7}[V_1]_{q}$}
\begin{center}
\begin{tabular}{c|ccc}
$\F_0=\R$ &  $\F_1=\R$ & $\F_1=\C$ & $\F_1=\HH$ \\
\hline
$q=0$
& $[V_0 \otimes_{\R} V_1]_{7}$ 
& $[V_0 \otimes_{\R} rV_1]_{7}$
& $[V_0 \otimes_{\R}  rc'V_1]_{7}$ \\
$q=3$
& $0$
& $0$ 
& $[V_0 \otimes_{\R} V_1]_{2}$\\
$q=4$
& $0$ 
& $0$
& $[V_0 \otimes_{\R} V_1]_{3}$\\
$q=7$
& $[V_0 \otimes_{\R} V_1]_{6}$ 
& $0$
& $0$
\end{tabular}  
\end{center} 
\end{table}
\end{lemma}
\proof
We prove some in the tables and the others will be left to the reader.

Let $\F_0=\R$.
Since $\beta_{\R, 7}\beta_{\R, 7}=\beta_{\R, 8+6}$, 
$\beta_{\R,7}\beta_{\HH, 4}=\beta_{\HH, 8+3}$,
and $\beta_{\R,7}\beta_{\HH, 3}=\beta_{\HH, 8+2}$ by Remark \ref{relation},
if $\F_1=\R$ then we have $[V_0]_7[V_1]_7=[V_0 \otimes _{\R} V_1]_6$, 
and if $\F_1=\HH$ then we have  $[V_0]_7[V_1]_4=[V_0 \otimes _{\R} V_1]_3$, 
$[V_0]_7[V_1]_3=[V_0 \otimes _{\R} V_1]_2$.

Let $\F_0=\HH$.
We fix the canonical isomorphism 
$Cl(\R^4)\otimes_{\R}Cl(\R^4)\cong Cl(\R^8)$ as $\Z/2$-graded real algebras
and consider $\Delta_4$ as a $\Z/2$-graded quaternionic $Cl(\R^4)$-module.
Then we have 
$$\Delta_4 \otimes _{\R} \Delta_4\cong \HH \otimes _{\R}\Delta_8 $$
as $\Z/2$-graded $\HH \otimes _{\R} Cl(\R^8)\otimes_{\R}\HH$-modules,
where $\HH$ in the right hand side is regarded
as the quaternionic bimodule by multiplication
from the both sides. 
Thus if $\F_1=\HH$, we have 
$[V_0]_3[V_1]_4=[V_0 \otimes _{\HH} V_1]_7$
and $[V_0]_3[V_1]_3=[V_0 \otimes _{\HH} V_1]_6$
since
$\beta_{\R, 7}\beta_{\R, 7}=\beta_{\R,8+6}$.
\qed

\section{Real and symplectic $K$-theory}\label{real-symplectic}

In this section we briefly review Real and symplectic $K$-theory.
We also discuss natural maps between other $K$-theories,
which will be used to construct spin $\pin$-vector bundles later.

Let $B$ be a compact Real space, i.e. a compact space with involution.
By the definition  in \cite{Atiyah-R} a Real vector bundle $V$ is
a  complex vector bundle $V$ together with
an anti-linear lift of the involution with its square equal to $1$.
We let $KR(B)$ be the Grothendieck group of
Real vector bundles over $B$.
According to Dupont \cite{Dupont} a symplectic vector bundle $V$ 
is defined to be a complex vector bundle $V$ together with
an anti-linear lift of the involution with its square equal to $-1$.
We let $Ksp(B)$ be the Grothendieck group of
symplectic vector bundles over $B$,
so that $KM(B)=KR(B)\oplus Ksp(B)$ 
becomes a $\Z/2$-graded ring in the obvious way.
These definitions then extend to a locally compact Real space $B$
in the usual way.
In the following we denote by $-1$ the involution on $B$
and regard $B$ as a $\{\pm 1\}$-space.

\begin{enumerate}
\item 
Let $V$ be a quaternionic $\{ \pm 1\}$-vector bundle over $B$.
Then the simultaneous action of $-1$ and $j \in \HH$ 
on $V$ defines a symplectic structure on $V$.
We denote by $\sigma V$ this symplectic vector bundle.
Thus we obtain a natural homomorphism 
$\sigma: KSp_{\{\pm 1\}}(B) \to Ksp(B)$.

\item
Let $V$ be a symplectic vector bundle $V$ over $B$.
From the symplectic action and complex scalar multiplication
on $V$ we obtain a real $\pin$-action on $V$.
We denote by $\vartheta V$ this real $\pin$-vector bundle.
Then we have a natural homomorphism 
$\vartheta: Ksp(B) \to KO_{\pin}(B)$.
Let $\tau: \pin \to \{\pm 1 \} \times Sp(1)$ 
be the inclusion defined by $\tau(t)=(1,t)$ and $\tau(j)=(-1, j)$.
Then it is obvious $\tau^*\rc V=\vartheta \sigma V$ for
a quaternionic $\{ \pm 1\}$-vector bundle $V$.
\item
Let $V$ be a symplectic vector bundle over $B$.
If we ignore the symplectic structure on $V$,
we get a complex vector bundle over $B$, which we write as $\kappa V$.
Then we obtain a natural homomorphism 
$\kappa: Ksp(B) \to K(B)$.
\end{enumerate}

\begin{remark}
The groups $KSp(B)$ and $Ksp(B)$ are similar in notation. 
However the latter is defined only for a Real space $B$,
and when the involution on $B$ is trivial
$KSp(B)$ is naturally isomorphic to $Ksp(B)$.
\end{remark}

Let $\R^{p,q}=\R^q \oplus \tilde \R^p$. 
We suppose $p \equiv q \pmod 8 $
and denote by 
$\beta_{\C}(\R^{p,q}) \in KR(\R^{p,q})$ 
the Bott class defined in \cite{Atiyah}.
Then the argument in \cite[Theorem 6.3]{Atiyah}
shows that multiplication by $\beta_{\C}(\R^{p,q})$ 
induces an isomorphism
$KM(B) \to KM(B \times \R^{p,q})$.

Now we recall some results in \cite{Dupont}.
Let $\beta_{\HH}(\tilde \R^{8k+4}) \in KSp_{\{\pm 1 \}}(\tilde \R^{8k+4})$
be the Bott class,
so $\sigma\beta_{\HH}(\tilde \R^{8k+4}) \in Ksp(\tilde \R^{8k+4})$.
Then by the uniqueness (up to grading) of the irreducible $\Z/2$-graded module
for the Real Clifford algebra 
associated to $\R^{8k,0}=\tilde \R^{8k}$ (see \cite[Section 4]{Atiyah-R}),
we see the square $(\sigma\beta_{\HH}(\tilde \R^{8k+4}))^2$
coincides with $\beta_{\C}(\R^{16k+8,0})$.
Thus multiplication by $\sigma\beta_{\HH}(\tilde \R^{8k+4})$
induces isomorphisms
$$KR(B) \to Ksp(B \times \tilde \R^{8k+4}),
\quad Ksp(B) \to KR(B \times \tilde \R^{8k+4}).$$
Combining them with the periodicity isomorphisms
\begin{align*}
KR(B \times \tilde \R^q)
\cong  KR(B \times \R^{8m,8m-q})
\cong  KR(B \times \R^{8m-q})
\cong  KO^q(B),
\end{align*}
we have the following lemma \cite{Dupont}:
\begin{lemma}\label{forgetful1}
The group $Ksp(\tilde\R^{n})$ is non-trivial only in the cases
\begin{align*}
Ksp(\tilde \R^{8k})\cong Ksp(\tilde\R^{8k+4}) \cong \Z, \quad
Ksp(\tilde \R^{8k+2})\cong Ksp(\tilde \R^{8k+3})\cong \Z/2\Z,
\end{align*}
and they are respectively generated by
\begin{align*}
[\HH] \beta_{\C}(\R^{8k, 0}), \,
\sigma\beta_{\HH}(\tilde \R^{8k+4}), \,
i_{8k+2}^*\sigma\beta_{\HH}(\tilde \R^{8k+4}),
\,
i_{8k+3}^*\sigma\beta_{\HH}(\tilde \R^{8k+4}),
\end{align*}
where $i_q: \tilde \R^q \times \{0\} \to \tilde \R^{8k+4}$ 
\rm{($q=8k+2, 8k+3$)} is the inclusion.
In addition, we have 
$$\kappa\sigma\beta_{\HH}(\tilde \R^{8k+4})=\beta_{\C,8k+4}.$$
\end{lemma}

\section{Irreducible representations of $\dpin$}\label{irreducible}

In this section we list up some properties of
real irreducible representations of the group $\dpin$
obtained as a central extension of $\pin$ by the group of order $2$.
\begin{definition}
Let $C_4=\langle j\rangle$ be the cyclic group of order $4$ generated by
$j$. We define an action of $C_4$ on $U(1)$
by $j(t)=t^{-1}$ and denote by $\Groupp$ the semi-product
$C_4 \ltimes U(1)$.
Then the quotient group $\Groupp/\{1,(j^2, -1)\}$ is isomorphic to $\pin$.
\end{definition}

We have the exact sequence
$$
1 \to U(1) \to \dpin \to C_4 \to 1
$$
and the center of $\dpin$ consists of the four elements
$\{\pm 1, (j^2, \pm 1)\}$.
We write $t=(1,t)$, $j^p=(j^p,1)$ below.

\begin{lemma}
The real irreducible representations of $\dpin$
are classified as follows:
\begin{enumerate}
\item
$\R$: the trivial $1$-dimensional representation,
on which we have
$$
t=j=\id.
$$
\item
$\tilde{\R}$: the non-trivial $1$-dimensional
representation, on which we have
$$
t=\id,\quad j=-\id.
$$
\item
$\C_0$: the $2$-dimensional irreducible representation of
$C_4$, on which we have
$$
t=\id ,\quad j=\sqrt{-1} \id,
$$
where we identify $\C$ with the complex numbers.
\item
$\D_m$ $(m\geq 1)$: the $2$-dimensional
 representation that satisfies
$$
tr(t|\D_m)=t^m+ t^{-m},\quad tr(j|\D_m)=0,\quad j^2=\id.
$$
\item
$\HH_m$ $(m\geq 1)$: the $4$-dimensional
 representation that satisfies
$$
tr(t|\HH_m)=2(t^m+ t^{-m}),\quad tr(j|\HH_m)=0,\quad j^2=-\id.
$$
\end{enumerate}
\end{lemma}
For $m=0$ we define $\D_0$ and $\HH_0$ to be
$$
\D_0:=\R \oplus \tilde{\R},\quad \HH_0:=\C_0 \oplus \C_0
$$
so that the defining properties of $\C_m$ and $\HH_m$
are still satisfied.

\begin{lemma}\label{repofPin^-(2)}
\begin{enumerate}
\item
The real irreducible representations of $\pin$ are
$\R, \tilde{\R}, \D_{2n}$ and $\HH_{2n-1}$ for $(n\geq 1)$.
\item
The real irreducible representations of $C_4$ are
$\R, \tilde{\R}$ and $\C_0$.
\end{enumerate}
\end{lemma}
\proof
The latter part is obvious.  The former part is a consequence of
the following observation.
The irreducible representation of $\dpin$ is that of $\pin$
if and only if the action of $j^2$ coincides with
that of $-1 \in U(1)$.
\qed

For an irreducible representation space $V$,
the ring $\End^{\dpin}(V)$ of the $\dpin$-invariant endomorphisms
is a field. We have the following three cases.
\begin{itemize}
\item
Real case: $\End^{\dpin}(V) \cong \R$ for
$V=\R, \tilde{\R}$ and $\D_k$ $(k\geq 1)$.
We denote by $\Irr_{\R}$  the set consisting of
these representations.
\item
Complex case: $\End^{\dpin}(V) \cong \C$ for
$V=\C_0$. We denote by $\Irr_{\C}$  the set consisting of
this representation.
\item
Quaternionic case:
$\End^{\dpin}(V) \cong \HH$ for
$V=\HH_k$ $(k\geq 1)$.
We denote by $\Irr_{\HH}$  the set consisting of
these representations.
\end{itemize}

To obtain the decomposition (\ref{KOG}) for $\dpin$
we need to fix an isomorphism $\End^{\dpin}(V) \cong \C$ 
for each $V \in \Irr_{\C}$.
In our case the only element in $\Irr_{\C}$ is $\C_0$.
To fix an isomorphism $\End^{\dpin}(\C_0) \cong \C$ we use the identification 
of $\C_0$ with the complex numbers in the definition of $\C_0$.
Then we have

\begin{corollary}\label{KO2}
There are canonical isomorphisms
\begin{align*}
KO_{\dpin}^0(\pt) & 
\cong\Z [\R] \oplus \Z[\tilde \R] \oplus
\bigoplus_{m=1}^{\infty}\Z[\D_m] \oplus \Z [\C_0] \oplus
\bigoplus_{n=1}^\infty \Z [\HH_n],\\
KO_{\dpin}^2(\pt) & \cong \Z [\C_0]_2 \oplus 
\bigoplus_{m=1}^\infty \Z/2 [\HH_m]_2,\\
KO_{\dpin}^4(\pt) & 
\cong\Z [\R]_4 \oplus \Z[\tilde \R]_4 \oplus
\bigoplus_{m=1}^{\infty}\Z[\D_m]_4 \oplus  \Z [\C_0]_4 \oplus 
\bigoplus_{n=1}^\infty \Z [\HH_n]_4,\\
KO_{\dpin}^{6}(\pt)& \cong\Z/2[\R]_6 \oplus \Z/2[\tilde \R]_6 \oplus
\bigoplus_{m=1}^{\infty}\Z/2[\D_m]_6 \oplus \Z[\C_0]_6.
\end{align*}
\end{corollary}

The following first two tables describe the tensor product
of the irreducible representations over $\R$ or $\C$,
from which we use  Lemma \ref{decomposed}, \ref{product}
to get the table of the product of elements in $KO_{\dpin}^*(\pt)$ below.

\begin{table}[h]\label{table-spin1}
\caption{Tensor products over $\R$}
\begin{center}
\begin{minipage}[t]{11cm}
\hspace{-5mm}
\begin{minipage}[t]{5cm}
\begin{tabular}[t]{|cc|c|}
$V_0$ &  $V_1$ & $V_0 \otimes_{\R} V_1$ \\
\hline 
$\tilde{\R}$ &  $\tilde{\R}$ & $\R$  \\
$\tilde{\R}$ &  $r\C_0$      & $r\C_0$ \\
$\tilde{\R}$ & $\D_m$        & $\D_m$ \\
$\tilde{\R}$ & $rc'\HH_m$    & $rc'\HH_m$ \\
$r\C_0$      & $r\C_0$       & $2(\R \oplus \tilde{\R})$ \\
\end{tabular}  
\end{minipage}
\begin{minipage}[t]{5cm}
\begin{tabular}[t]{|cc|c|}
$r\C_0$      & $\D_m$        & $rc'\HH_m$ \\
$r\C_0$      & $rc'\HH_m$    & $4\D_m$ \\
$\D_m$       & $\D_n$        & $\D_{|n-m|} \oplus \D_{n+m}$ \\
$\D_m$       & $rc'\HH_n$    & $rc'(\HH_{|n-m|} \oplus \HH_{n+m})$\\
$rc'\HH_m$   & $rc'\HH_n$    & $4(\D_{|n-m|} \oplus \D_{n+m})$
\end{tabular}  
\end{minipage}
\end{minipage}
\end{center} 
\end{table}
\begin{table}[h]\label{table-spin2}
\caption{Tensor products over $\C$}
\begin{center}
\begin{minipage}[t]{11cm}
\begin{minipage}[t]{5cm}
\begin{tabular}[t]{|cc|c|}
$V_0$ &  $V_1$ & $V_0 \otimes_{\C} V_1$ \\
\hline 
$c\tilde{\R}$ & $c\tilde{\R}$ & $c\R$ \\
$c\tilde \R$ & $\C_0$         &  $t\C_0$ \\
$c\tilde \R$ & $t\C_0$        & $\C_0$ \\
$c\tilde{\R}$ & $c\D_m$       & $c\D_m$ \\
$c\tilde{\R}$ & $c'\HH_m$     & $c'\HH_m$  \\
$\C_0$      &  $\C_0$          & $c\tilde{\R}$ \\
$t\C_0$     & $\C_0$          & $c\R$ \\
\end{tabular}  
\end{minipage}
\begin{minipage}[t]{5cm}
\begin{tabular}[t]{|cc|c|}
$\C_0$      & $c\D_m$         &    $c'\HH_m$ \\
$t\C_0$     & $c\D_m$         & $c'\HH_m$\\
$\C_0$      & $c'\HH_m$       & $2c\D_m$\\
$t\C_0$   &   $c'\HH_m$       & $2c\D_m$ \\
$c\D_m$   &  $c\D_n$          & $c(\D_{|n-m|} \oplus \D_{n+m})$ \\
$c\D_m$  &   $c'\HH_n$        & $c'(\HH_{|n-m|} \oplus \HH_{n+m})$ \\
$c'\HH_m$ &  $c'\HH_n$        & $c(\D_{|n-m|} \oplus \D_{n+m})$
\end{tabular}  
\end{minipage}
\end{minipage}
\end{center} 
\end{table}

\begin{table}[h]\label{tensorKO2}
\caption{Table of $[V_0]_p[V_1]_q$}
\begin{center}
\begin{minipage}[t]{11cm}
\hspace{-5mm}
\begin{minipage}[t]{5cm}
\begin{tabular}[t]{|cc|c|}
$[V_0]_p$ & $[V_1]_q$ & $[V_0]_p[V_1]_q$ \\
\hline 
$[\R]$         & $[\C_0]_2$ & $[\C_0]_2$ \\           
$[\tilde{\R}]$ & $[\C_0]_2$ & $-[\C_0]_2$ \\ 
$[\C_0]$       & $[\C_0]_2$  & $0$ \\           
$[\D_n]$       & $[\C_0]_2$ & $[\HH_n]_2$ \\ 
$[\HH_n]$  &     $[\C_0]_2$  & $0$ \\         
$[\R]$      &     $[\HH_m]_2$ & $[\HH_m]_2$ \\ 
$[\tilde{\R}]$ &   $[\HH_m]_2$ & $[\HH_m]_2$\\
$[\C_0]$    &    $[\HH_m]_2$   & $0$ \\
$[\D_n]$      &   $[\HH_m]_2$   & $[\HH_{|n-m|}]_2+[\HH_{n+m}]_2$ \\
$[\C_0]_2$   &    $[\C_0]_2$   & $[\tilde \R]_4 -[\R]_4$ \\
$[\R]_4$    &     $[\HH_m]_4$ &   $2[\HH_m]$\\
$[\tilde \R]_4$ & $[\HH_m]_4$ & $2[\HH_m]$ \\
$[\HH_m]_4$ &   $[\HH_n]_4$  &   $4([\D_{|n-m|}]+[\D_{n+m}])$\\
\end{tabular}  
\end{minipage}
\qquad\qquad\qquad
\begin{minipage}[t]{5cm}
\begin{tabular}[t]{|cc|c|}
$[\R]_4$   &   $[\R]_4$    &   $4[\R]$\\
$[\tilde \R]_4$ & $[\tilde \R]_4$ & $4[\R]$\\
$[\R]_4$ &     $[\tilde \R]_4$  & $4[\tilde \R]$\\
$[\R]_4$ &  $[\C_0]_2$ & $2[\C_0]_6$ \\
$[\tilde{\R}]_4$ & $[\C_0]_2$ & $-2[\C_0]_6$ \\
$[\C_0]_4$ & $[\C_0]_2$  & $0$ \\
$[\D_n]_4$ & $[\C_0]_2$ & $0$ \\
$[\HH_n]_4$ & $[\C_0]_2$ & $0$ \\
$[\R]_4$  & $[\C_0]_6$   & $2[\C_0]_2$\\
$[\tilde{\R}]_4$ & $[\C_0]_6$ & $-2[\C_0]_2$\\
$[\C_0]_4$ & $[\C_0]_6$ &  $0$\\
$[\D_n]_4$ & $[\C_0]_6$ & $0$\\
$[\HH_n]_4$ & $[\C_0]_6$ & $0$
\end{tabular}  
\end{minipage}
\end{minipage}
\end{center}
\end{table}

\section{Calculations of the Euler classes}\label{calculations}

Let $\C_{(i)}$ be the complex representation of $C_4$
that satisfies $tr(j \vert \C_{(i)})=tr(j \vert \C_1)^i$.
Then the irreducible complex representation of $C_4$ 
consists of $\C_{(0)}, \C_{(1)}, \C_{(2)}$ and $\C_{(3)}$.

\begin{definition}
We put a spin $C_4$-structure on $\tilde \R^{2m}=\tilde \R^{\oplus 2m}$
and a spin $\pin$-structure on $\HH_1^m=\HH_1^{\oplus m}$
in the following way\rm{:}
\begin{itemize}
\item The representation $\C_{(1)}$ 
gives a square root of the action $C_4$ on $\C_{(2)}$,
which induces a spin $C_4$-structure on $r\C_{(2)}=\tilde \R^2$.
\item The $\pin$-action on $\HH_1$ factors through
the homomorphism $\pin \subset Sp(1)=\Spin(3) \to SO(4)$.
Hence the lift to $\Spin(4)$ gives a spin $\pin$-structure on $\HH_1$.
\item We put a spin $C_4$-structure on $\tilde \R^{2m}$
and a spin $\pin$-structure on $\HH_1^m$ 
as the direct sums of $\tilde \R^2$ and $\HH_1$ respectively.
\end{itemize}
\end{definition}
Then by the formulae (\ref{Bott-c}), (\ref{Euler-c})
we have
\begin{align*}
c(\beta(\tilde{\R}^2))
&=[\C_{(3)}]\beta_{\C}(\C_{(2)}),\quad
c(e(\tilde{\R}^2))=[\C_{(3)}] -[\C_{(1)}].
\end{align*}
Since $KO_{C_4}^{2}(\pt)\cong \Z[\C_0]_2$,
$K_{C_4}^{2}(\pt) \cong R(C_4)$
and $c[\C_0]_2=[\C_{(1)}] -[\C_{(3)}]$ by Lemma \ref{complexi},
we use Table 9 to obtain
\begin{lemma}\label{eR}
Consider the above spin $C_4$-structure on $\tilde \R^{2m}$.
Then we have 
$c(\beta(\tilde \R^{2m})) =[\C_{(3m)}]\beta_{\C}(\C_{(2)}^{\oplus m})$,
$c(e(\tilde \R^{2m})) =([\C_{(3)}]-[\C_{(1)}])^m$ and
\begin{align*}
e(\tilde \R^{2m}) = & 
\begin{cases}
(-1)^{m/2}2^{m-1}([\R]_{2m}-[\tilde \R ]_{2m}) & \text{$m \equiv 0 \mod 4$},\\
(-1)^{(m+1)/2}2^{m-1}[\C_0]_{2m} & \text{$m \equiv 1 \mod 4$},\\
(-1)^{m/2}2^{m-2}([\R]_{2m}-[\tilde \R ]_{2m}) & \text{$m \equiv 2 \mod 4$},\\
(-1)^{(m+1)/2}2^{m-1}[\C_0]_{2m} & \text{$m \equiv 3 \mod 4$}.
\end{cases}
\end{align*}
\end{lemma}
On the other hand the complexification of $e(\HH_1)$ is
\begin{align*}
c(e(\HH_1))&=[\Lambda^{\even}(c'\HH_1)]-[\Lambda^{\odd}(c'\HH_1)]
=2[\C_{(0)}] -[c'\HH_1],
\end{align*}
by (\ref{Euler-c})
and $c: KO_{\pin}^{4}(\pt) \to K_{\pin}^{4}(\pt)$ is injective
by Lemma \ref{complexi}, so we obtain
\begin{lemma}\label{eH}
Consider the above spin $\pin$-structure on $\tilde \HH^{m}$.
Then we have
$c(e(\HH_1))=2[\C_{(0)}] -[c'\HH_1]$ and
\begin{align*}
e(\HH_1)  & = [\R]_4 -[\HH_1]_4,\\
e(\HH_1)^{2m} & = a_m [\R] + b_m [\tilde \R] +
\sum_{1 \le t < 2m, \even }{4m \choose 2m-t}[\D_{t}] \\
 & -\dfrac 12\sum_{1 \le t < 2m, \odd}{4m \choose 2m-t}[\HH_{t}]
+ [\D_{2m}],
\end{align*}
where the integers $a_m$, $b_m$
satisfy $a_m-b_m=2^{2m}$, $a_m+b_m={4m \choose 2m}$.
\end{lemma}

Let $\WH^+(k)$ be the quaternionic $\{\pm 1\}$-vector bundle
over the $\{\pm 1\}$-sphere
$\tilde \R^{8k+4} \cup \{ \infty\}$
which is the extension of the product quaternionic 
$\{\pm 1\}$-vector bundle
$\tilde \R^{8k+4} \times \Delta^{+}_{8k+4}$
obtained by using the Clifford multiplication $c$
as the trivialization  of 
$\tilde \R^{8k+4} \times \Delta^{+}_{8k+4}$
over $\tilde \R^{8k+4} \setminus \{0\}$.
If we exchange $\Delta^+_{8k+4}$ for $\Delta^-_{8k+4}$
and replace $c$ by $c^{-1}$
then we obtain another quaternionic $\{\pm 1\}$-vector bundle $\WH^-(k)$.
These bundles then satisfy
\begin{align*}
[\WH^+(k)] & =[\Delta^-_{8k+4}]
+\beta_{\HH}(\tilde \R^{8k+4}), \quad
[\WH^-(k)] =[\Delta^+_{8k+4}]
-\beta_{\HH}(\tilde \R^{8k+4})
\end{align*}
as elements in $KSp_{\{\pm 1\}}(\tilde \R^{8k+4} \cup \{ \infty\})$
(see Section \ref{some-properties}).
As shown in Section \ref{some-properties},
we have a natural spin $\{ \pm 1\} \times Sp(1)$-structure
on $\rc \WH^{\pm}(k)$,
so we have a spin $\pin$-structure on
$\vartheta\sigma\WH^{\pm}(k)=\tau^*\rc\WH^{\pm}(k)$.
In the following 
the Euler class 
$e(\vartheta\sigma\WH^{\pm}(k))
 \in KO_{\pin}^{2^{4k+2}}(\tilde \R^{8k+4} \cup \{\infty
\})$
will be also denoted 
$e(\WH^{\pm}(k))$ for simplicity.
\begin{remark}
In the above construction
the group $\pin$ was defined as the subgroup of $\{\pm 1\} \times Sp(1)$
generated by $\{1\} \times U(1)$ and $(-1,j)$.
However we also used the symbol $j \in \pin$ instead of $(-1,j)$.
To avoid confusion in the rest of this section
we denote by $j_{\HH}$ the element $j \in \HH$,
while we denote by $j_4$ the element $(-1,j_{\HH}) \in \pin$.
Thus $C_4$ is the subgroup of $\pin$ generated by $j_4$.
\end{remark}
\begin{lemma}\label{E-H}
Under the canonical decomposition 
$KO_{\pin}^{2^{4k+2}}(\tilde \R^{8k+4} \cup \{ \infty\})
\cong KO_{\pin}^{2^{4k+2}}(\{ \infty\})
\oplus KO_{\pin}^{2^{4k+2}}(\tilde \R^{8k+4})$
we have
\begin{align*}
e(\WH^{+}(k))
& =(-[\tilde \R]e(\HH_1))^{2^{4k}}+\gamma_k\beta(\tilde \R^{8k+4}), \\
e(\WH^{-}(k))
& =(-e(\HH_1))^{2^{4k}}-\gamma_k\beta(\tilde \R^{8k+4}),
\end{align*}
where $\gamma_0=-[\tilde \R]$,
and $\gamma_k \in KO_{\pin}^4(\pt)$ \rm{($k \ge 1$)} 
is in the subgroup generated by
$[\R]_4+[\tilde \R]_4$, $[\C_0]_4$,
$[\D_{2m}]_4$ and $[\HH_{2m-1}]_4$ \rm{($m \ge 1$)}.
\end{lemma}

\proof
We first consider the action of the subgroup $C_4$ in $\pin$.
The action of $j_{\HH} \in \HH$ then induces a complex $\{ \pm 1\}$-structure 
on the quaternionic $\{\pm 1\}$-structure on $\WH^+(k)$,
so we denote by $c''\WH^+(k)$ this complex $\{ \pm 1\}$-vector bundle.
Let $i_{C_4}: C_4 \to \{ \pm 1\} \times U(1)$ be the inclusion
derived from the inclusion $\tau: \pin \to \{ \pm 1\} \times Sp(1)$.
Since the restriction $c'\Delta_{8k+4}$ of the scalars to
$\C=\R \oplus j\R$ 
is isomorphic to $\Lambda^*(\C^{4k+2})$
as $\Z/2$-graded complex $Cl(\R^{8k+4})$-modules,
we see $i_{C_4}^*bc''\WH^+(k)$ is isomorphic to 
$\C_{(1)}\otimes \hat V^+$ as complex $C_4$-vector bundles,
where $V=\C_{(2)}^{\oplus 4k+2}$
and $\hat V^+$ is the complex $C_4$-vector bundle
defined in Section \ref{some-properties}.
If we denote by $j_{C_4}: C_4 \to \pin$ the inclusion,
then from Proposition \ref{Euler-class-0}, \ref{Euler-class-1}
we obtain
\begin{align*}
& j_{C_4}^*c(e(\vartheta\sigma\WH^+(k)))\\
& =[\C_{(3)}]^{2^{4k}}e_{\C}(i_{C_4}^*bc''\WH^+(k))\\
& =[\C_{(3)}]^{2^{4k}}\left( e_{\C}(\C_{(3)})^{2^{4k+1}}
+\mu_{4k+2}(\C_{(1)}, \C_{(2)}, \cdots, \C_{(2)})
\beta_{\C}(\C_{(2)}^{\oplus 4k+2})\right)\\
& =\left(-[\C_{(2)}]e_{\C}(\C_{(1)}\oplus \C_{(3)})\right)^{2^{4k}}
+\gamma'_k\beta_{\C}(\C_{(2)}^{\oplus 4k+2}),
\end{align*}
where $\gamma'_0 =-1$ and $\gamma'_k \in R(C_4)$
is seen to satisfy $tr(j_4 | \gamma'_k)=0$ ($k\ge 1$).

We next consider the action of the subgroup $U(1)$.
Let $i_{U(1)}: U(1) \to \{\pm 1 \} \times U(1)$ be the inclusion.
Then $i_{U(1)}^*bc'\WH^+(k)$ is isomorphic to $\C_1(t) \otimes \hat V^+$
as complex $U(1)$-vector bundles,
where $V=\C_0(t)^{\oplus 4k+2}$.
If we denote by $j_{U(1)}: U(1) \to \pin$ the inclusion,
we use Proposition \ref{Euler-class-0} to get
\begin{align*}
j_{U(1)}^*c(e(\vartheta\sigma\WH^+(k)))
&= [\C_{-1}(t)]^{2^{4k}}e_{\C}(i_{U(1)}^*bc'\WH^+(k))\\
&=\left(-e_{\C}(\C_1(t) \oplus \C_{-1}(t))\right)^{2^{4k}}
+\gamma''_k\beta_{\C}(\C_0(t)^{\oplus 4k+2}),
\end{align*}
where $\gamma''_0=-1$ and $\gamma''_k \in KO_{U(1)}^{4}(\pt)$.

Thus we conclude $\gamma_0=-[\tilde \R]$ by Lemma \ref{eR}.
When $k \ge 1$,
Lemma \ref{complexi} can be applied
to show that the multiplicity of $[\R]_4$ in $\gamma_k$
should be equal to that of $[\tilde \R]_4$,
since $\Tr(j_4|c(\gamma_n))=\Tr(j_4 | \gamma'_n)=0$.
The calculation of the Euler class $e(\WH^-(k))$ is similar.
\qed

If we apply Proposition \ref{Euler-class-2} in the above proof,
we can compute the quantity $\gamma_k \in KO_{\pin}^4(\pt)$ ($k \ge 1$), 
thus the Euler classes $e(\hat \HH^{\pm}(k))$ inductively.
As mentioned in Section \ref{some-properties} we do not carry out it.
However we need the following corollary
obtained from Lemma \ref{E-H} together with Table 9 
and Lemma \ref{eR}.
\begin{corollary}\label{keyrelation}
\begin{align*}
e(\tilde \R^2)e(\WH^{\pm}(0))^2 &=e(\tilde \R^2)e(\HH_1)^2
\in KO_{\dpin}^2(\tilde \R^{8k+4} \cup \{ \infty\}),\\
e(\tilde \R^2)e(\WH^{\pm}(k))
& = e(\tilde \R^2)e(\HH_1)^{2^{4k}}
\in KO_{\dpin}^2 (\tilde \R^{8k+4} \cup \{ \infty\})
\quad \text{\rm{($k \ge 1$)}}.\\
\end{align*}
\end{corollary}

\section{$KO_{\dpin}^*(\tT^n)$}\label{KO-torus}
Let $A$ be a $\dpin$-space.
Suppose the $\dpin$-action is not free.
We take a fixed point $a_0$ in $A$.
Let $S$ be a subset of $\{1, \cdots, n \}$.
Then the $\dpin$-space $A^S= \Map(S,A)$ can be embedded
$\dpin$-equivariantly into $A^n$ by the map
$$h_S : A^{S}\overset{\cong}\longrightarrow
\{ (x_1, \cdots , x_n) \in A^n
\mid x_i =x_0 \, (i \not\in  S) \} \subset A^n.$$
Let $\pi_S: A^n \to A^S$ be the projection.
For a $\dpin$-invariant open neighborhood $U$ around $a_0$ in $A$
we let $i_S: U^S \to A^S$ be the corresponding inclusion.
For the subset $S=\{1, \cdots, m\}$ ($m \le n$)
we write $A^m=A^S$, $\pi=\pi_S: A^n \to A^m$
and $h=h_S: A^m \to A^n$ for simplicity.

We recall that $\pin$ acts on $\tilde \R$ via 
the map $\pin \to \{\pm 1\}$,
hence $\dpin$ acts on $\tilde \R$ via the projection
$\dpin \to \pin$, 
by which the torus $\tT^n =(\tilde \R/\Z)^n$ becomes a $\dpin$-space.
We now take a fixed point $t_0$ in $\tT^1$ 
and consider the following exact sequence
\begin{align*}
\to KO_{\dpin}^q((\tT^{n}, \tT^{n-1})\times \tilde \R^{p})
\overset{j^*}{\to}KO_{\dpin}^q(\tT^{n}\times \tilde \R^{p})
\overset{h^*}{\to}KO_{\dpin}^{q}(\tT^{n-1}\times \tilde \R^{p})
\to,
\end{align*}
where the first term
is $KO_{\dpin}^q((\tT^{n}, \tT^{n-1})
\times (\tilde \R^{p} \cup \{\infty\}, \{\infty\}))$,
which is identified with
$KO_{\dpin}^q(\tT^{n-1}\times \tilde \R^{1+p})$
by excision.
Then $j^*$ is identified with the push-forward map
$i_!: KO_{\dpin}^q(\tT^{n-1}\times \tilde \R^{1+p})
\to KO_{\dpin}^q(\tT^n\times \tilde \R^{p})$
induced from a $\dpin$-open embedding $i: \tilde \R \to \tT^1$
onto a neighborhood of $t_0$.
Since $\pi^*$ give a right-inverse of $h^*$, 
the sequence is split and we obtain an isomorphism
\begin{equation*}
i_! + \pi^*:
KO^q_{\dpin}(\tT^{n-1}\times \tilde \R^{1+p})
\oplus
KO^q_{\dpin}(\tT^{n-1}\times \tilde \R^{p})
\rightarrow KO_{\dpin}^q(\tT^n \times \tilde \R^p).
\end{equation*}
By induction on the cardinal number $|S|$ of $S$ we see
\begin{lemma}
The following map is an isomorphism\rm{:}
\begin{equation}\label{torus}
\sum_{S \subset \{1, \cdots , n \}}\pi_S^*(i_S)_!:
\bigoplus_{S}KO^q_{\dpin}(\tilde \R^S\times \tilde \R^p)
\rightarrow KO_{\dpin}^q(\tT^n \times \tilde \R^p ).
\end{equation}
\end{lemma} 
\begin{remark}
It is easy to see
the decomposition of (\ref{torus}) still holds
when we replace $KO_{\dpin}$ by other $K$-groups as
$KSp_{\{\pm 1\}}$, $Ksp$, $KR$ and $K$.
We shall also use these isomorphisms later.
\end{remark}
For each subset $S \subset \{1, \cdots, n \}$
we fix an identification $S$ with $\{1, \cdots, |S|\}$ as sets.
Then it induces an identification $\tilde \R^S$ with $\tilde \R^{|S|}$.
If $|S|$ is even,
we define the Bott class
$\beta(\tilde \R^{S}) \in  KO_{C_4}^{|S|}(\tilde \R^S)
\subset KO_{\dpin}^{|S|}(\tilde \R^S)$
to be $\beta(\tilde \R^{|S|})$
under this identification.
The Bott periodicity theorem then asserts
that $KO_{\dpin}^q(\tilde \R^S)$ is freely generated
by the Bott class $\beta(\tilde \R^S)$
as a $KO_{\dpin}^{q-|S|}(\pt)$-algebra.

Let $p$ be a non-negative integer and suppose $|S|$ is even.
Then by definition we have
$$
e(\tilde \R^{p})\beta(\tilde \R^S)
=i^*\beta(\tilde \R^p \oplus \tilde \R^S)
\in KO_{\dpin}^{p+|S|}(\tilde \R^S)\subset KO_{\dpin}^{p+|S|}(\tT^n),
$$
where $i:\tilde \R^S \to \tilde \R^p \oplus \tilde \R^S $
is the inclusion.
More generally if $p +|S|$ is even
then the right hand side is still defined,
so we may write it as $e(\tilde \R^{p})\beta(\tilde \R^S)$.
Then one can immediately see the following product formula:
\begin{equation}\label{Bott-product}
e(\tilde \R^p)\beta(\tilde \R^{S})
e(\tilde \R^{p'})\beta(\tilde \R^{S'})
= e(\tilde \R^{p+p'+\vert S \cap S'\vert})
\beta(\tilde \R^{S \cup S'})
\end{equation}
for any subsets $S, S' \subset \{1, \cdots, n\}$
and integers $p, p' \ge 0$ with $|S|+p$, $|S'|+p'$ even.
\begin{remark}
Even if $p$ or $|S|$ is not even
it is possible to define $\beta(\tilde \R^S)$ and $e(\tilde \R^p
)$
separately 
so that the above two (equivariant) product formulae hold,
if one uses the notion of $KO$-group with local coefficients 
\cite{Donovan-Karoubi}.
However in this paper we do not introduce this notion,
since the above definition is sufficient in our calculation.
\end{remark}

\section{Proof of Theorem~\ref{main1}}\label{Proof-of-1}

Let $l$ be a positive even integer.
Let $V_0=\tT^n\times \tilde \R^{x}$, $W_0=\tT^n \times \tilde \R^{x+l}$
be the product bundles over $\tT^n$,
and let $V_1$, $W_1$ be symplectic bundles over $\tT^n$.
We assume the difference $[V_1]-[W_1] \in Ksp(\tT^n)$
satisfies the condition of Theorem~\ref{main1}.
Suppose we have a proper $\pin$-equivariant fiber-preserving map
$\varphi: V_0\oplus W_0 \to V_1 \oplus W_1$
which induces the identity on the base space $\tT^n$
and whose restriction $\varphi_0: V_0 \to W_0$ is given by 
the standard inclusion $\tilde \R^{x} \to \tilde \R^{x+l}$.

We use the projection $\dpin \to \pin$
to regard $V_1$, $W_1$ as $\dpin$-equivariant bundles.
Let $\dpin_0=U(1)$ be the subgroup in $\dpin$,
so the quotient $\dpin/\dpin_0$ is isomorphic to $C_4$.
Using the projections $\dpin \to C_4 \to \{\pm 1\}$
we regard $V_0$, $W_0$ as $C_4$-equivariant bundles,
and also $\dpin$-equivariant bundles.
In this setting $\varphi$ is a $\dpin$-equivariant map.

We first calculate $KO_{C_4}^*$-degree $\alpha_{\varphi_0}$
of $\varphi_0$. 
After stabilization 
by the identity map on $\tT^n \times  \tilde \R$ if necessary,
we may assume that $x$ is even.
Then $i^*\beta(\tilde \R^{x+l})
=\beta(\tilde \R^{x})e(\tilde \R^l)$,
and so we have 
$$\alpha_{\varphi_0}=e(\tilde \R^l) \in KO_{C_4}^l(\pt).$$

We next calculate the $KO_{\dpin}^*$-degree
$\alpha_{\varphi} \in KO_{\dpin}^{l-4k}(\tT^n)$.
In the decomposition of $Ksp(\tT^n)$
we take $\sigma\beta_{\HH}(\tilde \R^{S})$ 
as a generator of $Ksp(\tilde \R^{S}) \cong \Z$
when $\caF\equiv 4 \pmod 8$,
by which we identify $a_S$ with an integer.
On the other hand when $\caF\equiv 2,3  \pmod 8$
we choose an integer $a_S$ as an lift of $\overline a_S \in \Z/2$.
Then Lemma~\ref{forgetful1} implies that,
after stabilization 
by the identity map on some symplectic bundle over $\tT^n$,
we may assume that the symplectic bundles
$V_1$, $W_1$ take of the form:
\begin{align*}
V_1 &=\sigma\HH^{y+A+k},
\quad A=\sum _{\caF =8p+2, \, 8p+3, \, 8p+4,}|a_S|2^{4p},\\
W_1 & =\sigma\HH^y \oplus 
\bigoplus_{\caF =8p+4}
\pi_S^*f_{S}^*(\sigma
(\overbrace{
\WH^{\epsi(S)}(p)
 \oplus \cdots \oplus \WH^{\epsi(S)}(p)
}^{\text{$|a_S|$ times}})) \\
& \qquad\quad \oplus \bigoplus_{\caF =8p+2, \, 8p+3}
\pi_S^*f_{S}^*\tilde i_{|S|}^*(\sigma(\overbrace
{\WH^{\epsi(S)}(p)
\oplus \cdots \oplus \WH^{\epsi(S)}(p)}^{\text{$|a_S|$  times}})
\end{align*}
where $\tilde i_q$ ($q=8p+2, 8p+3$) is the extension of
the inclusion $i_q: \tilde \R^{q} \times \{ 0\} \to \tilde \R^{8p+4}$
to the one-point compactifications,
and where 
$f_{S}: T^{S} \to \tilde \R^{S} \cup \{\infty\}
\cong \tilde \R^{|S|} \cup \{\infty\}$
is a $\{\pm 1\}$-equivariant map obtained by shrinking the complement of 
a neighborhood of the fixed point $(t_0, \cdots, t_0)$ in $T^{S}$ to the
one point $\{ \infty \}$,
and where $\epsilon(S)=$''$+$'' when $a_S > 0$,
and $\epsilon(S)=$''$-$'' when $a_S < 0$.
By Corollary \ref{keyrelation}
the second equation of Lemma~\ref{keylemma} becomes
\begin{align}\label{keyequation}
\alpha_{\varphi} e(\HH_1)^{y+A+k}
& = e(\tilde \R^l)e(\HH_1)^{y+A-A_1}
\prod_{\caF=4, \,\, \overline a_S \ne 0}
\pi_S^*f_{S}^*e(\WH^{\epsi(S)}(p)) \notag\\
& \qquad \qquad \qquad \quad
\times \prod_{\caF=2,3, \,\, \overline a_S \ne 0}
\pi_S^*f_{S}^*\tilde i_{|S|}^*
e(\WH^{\epsi(S)}(p)),
\end{align}
where $A_1$ is the cardinal number of subsets $S$
satisfying $|S|=2,3, 4$ and $\overline a_S \ne 0$.

\begin{proposition}\label{calc-KO-degree}
Let $l$ be a positive even integer.
Suppose $\alpha_{\varphi} \in KO_{\dpin}^{l-4k}(\tT^n)$
satisfies the equation \rm{(\ref{keyequation})}.
Let $(\alpha_{\varphi})_S \in KO_{\dpin}^{l-4k}(\tilde \R^S)$ 
be the component of 
$\alpha_{\varphi} \in KO_{\dpin}^{l-4k}(\tT^n)$
in the decomposition of $KO_{\dpin}^{l-4k}(\tT^n)$.
If $\caS$ is even then we have
\begin{align*}
(\alpha_{\varphi})_S &=
\begin{cases}
\varepsilon
N_S 2^{l/2-k-|S|/2-1} ([\R]_{d}-[\tilde \R]_{d})
\beta(\tilde \R^{S}), &  d \equiv 0 \mod 8, \\
\varepsilon
N_S 2^{l/2-k-|S|/2-1} ([\C_0]_{d}+\text{(torsion)})
\beta(\tilde \R^{S}), &  d \equiv 2 \mod 8, \\
\varepsilon
N_S 2^{l/2-k-|S|/2-2}([\R]_{d}-[\tilde \R]_{d})
\beta(\tilde \R^{S}), &  d \equiv 4 \mod 8, \\
\varepsilon
N_S 2^{l/2-k-|S|/2-1} \left([\C_0]_{d}+ \text{(torsion)})\right)
\beta(\tilde \R^{S}), &  d \equiv 6 \mod 8,
\end{cases}
\end{align*}
where $d=l-4k-|S|$ and $\varepsilon = \pm 1$.
Moreover if $d \equiv 2 \pmod 8$ and $l \ge 4$ then 
$(\alpha_{\varphi})_S$ has no component in the torsion subgroup
and $N_S2^{l/2-k-\caS/2-2}$ is an integer.
\end{proposition}
To prove this we use the following lemma:
\begin{lemma}\label{divis}
Let $c \ge 0$ be even, $i \ge 0$ and $d$ even.
Suppose $\alpha \in KO_{\dpin}^{d+2i}(\tilde \R^{2i})$ satisfy the equation
\begin{align*}
\alpha e(\HH_1)^c
=\sum_{0 \le n < c+d/4}
a_n e(\HH_1)^n e(\tilde \R^{4c+d-4n})\beta(\tilde \R^{2i})
\in KO_{\dpin}^{4c+d+2i}(\tilde \R^{2i}),
\end{align*}
where $a_n \in \Z$ \rm{($0 \le n < c+d/4$)}.
Then we have
\begin{align*}
&\alpha = \\
&
\begin{cases}
\sum a_n (-1)^{c+d/4-n}2^{c+d/2-n-1}([\R]_{d}-[\tilde \R]_{d})
\beta(\tilde \R^{2i}), 
& d \equiv  0 \mod 8, \\
\sum a_n (-1)^{c+d/4-n+1/2}2^{c+d/2-n-1}
\left([\C_0]_{d} + \text{(torsion)} \right)
\beta(\tilde \R^{2i}), 
& d \equiv 2 \mod 8, \\
\sum a_n (-1)^{c+d/4-n}2^{c+d/2-n-2}([\R]_{d}-[\tilde \R]_{d})
\beta(\tilde \R^{2i}), 
& d \equiv 4 \mod 8, \\
\sum a_n (-1)^{c+d/4-n+1/2}2^{c+d/2-n-1}
\left( [\C_0]_{d}+ \text{(torsion)}\right)
\beta(\tilde \R^{2i}), 
& d \equiv 6 \mod 8.
\end{cases}
\end{align*}
Moreover if $d \equiv 2 \pmod 8$ 
and $a_n=0$ for all $n$ with $c +d/4 -1 < n$
then $\alpha$ has no component in the torsion subgroup
and $\sum a_n (-1)^{c+d/4-n+1/2}2^{c+d/2-n-2}$ is an integer.
\end{lemma}
\proof 
From Lemma \ref{eR} we see
the complexification $c(\alpha) \in K_{\dpin}^{d+2i}(\pt)\cong R(\dpin)$
satisfies
$$
\Tr(j^2|c(\alpha))=\Tr(t|c(\alpha))=0, \quad
\Tr(j|c(\alpha))=\sum a_n (-2i)^{2c+d/2-2n}2^{n-c},
$$
since $c(e(\HH_1))=2[\C_{(0)}]-[\HH_1]$ and
$c(e(\tilde \R^2))=[\C_{(3)}]-[\C_{(1)}]$.
This implies $c(\alpha)$ is in the subgroup $R(C_4)$ and
\begin{align*}
c(\alpha) 
=\sum a_n ([\C_{(3)}]-[\C_{(1)}])^{2c+d/2-2n}([\C_{(0)}]-[\C_{(2)}])^{n-c}.
\end{align*}
Since the kernel of 
$c: KO_{\dpin}^{d+2i}(\pt) \to K_{\dpin}^{d+2i}(\pt)$ is
the torsion subgroup by Lemma \ref{complexi},
we get the equation for $\alpha$.

When $d \equiv 2 \pmod 8$
Corollary~\ref{KO2} shows the torsion subgroup is generated by 
$[\HH_m]_{d}$ ($m \ge 1$).
Suppose $\alpha$ has a nontrivial $[\HH_m]_d$-component.
Let $m_0$ be the maximal of such $m$.
Note $[\HH_{m_0}]_{d}[\D_{c}]=[\HH_{|m_0-2c|}]_{d}+[\HH_{|m_0+2c|}]_{d}$.
It implies that $\alpha e(\HH^1)^c$ has a nontrivial 
$[\HH_{|m_0+2c|}]_{d}$-component.
On the other hand 
if we assume $a_n=0$ for all $n$ with $c +d/4 -1 < n$,
then we see from Lemma~\ref{eR}, \ref{eH} that
$$\sum_{ 0 \le n < c +d/4}
a_n e(\HH_1)^n e(\tilde \R^{4c+d-4n})
=\sum_n a_n(-1)^{c+d/4-n+1/2}2^{2c+d/2-n-1}[\C_0]_{4c+d},$$
which is a contradiction.

Moreover if $\sum a_n (-1)^{c+d/4-n+1/2}2^{c+d/2-n-1}$ were odd,
$$\alpha e(\HH_1)^{c}=\sum
a_n (-1)^{c+d/4-n+1/2}2^{2c+d/2-n-1}[\C_0]_{4c+d}+[\HH_{4c}]_{4c+d},
$$
which is also a contradiction.
\qed

{\it Proof of Proposition~\ref{calc-KO-degree}.}
We consider the right hand side of
(\ref{keyequation}) on $S \subset \{1, \cdots, n \}$.
Since $-[\tilde \R]e(\tilde \R^2)=e(\tilde \R^2)$,
we use (\ref{Bott-product}),
Lemma \ref{E-H} and Corollary~\ref{keyrelation} 
to get the expansion:
\begin{align*}
\varepsilon \sum_S \sum_{m \ge 0}N(S,m)
e(\HH_1)^{y+A-m}e(\tilde \R^{l+4m-\caS})\beta(\tilde \R^{S}),
\end{align*}
where $\varepsilon$ is $1$ or $-1$.
We may assume $y$ is large enough
and $y+A+k$ is even.
We then apply Lemma~\ref{divis} for
$c=y+A+k \ge 0, \, d=l-4k-|S|, \, 2i=|S|, \, n=y+A-m \ge 0, \,
a_n=N(S,m)$ and $\alpha=(\alpha_{\varphi})_S$.
Since $c+d/2-n =l/2-k-|S|/2+m$, $c+d/4-n =l/4-|S|/4+m \ge l/4 >0$,
and $c+d/4-n \ge 1$ for $l \ge 4$,
we have the conclusion.
\qed

Now we prove Theorem \ref{main1}.
Suppose the inequality were not satisfied.
Then by the standard inclusion 
$\tilde \R^{x+l} \to \tilde \R^{x+l+z}$ for some $z \ge 0$
we may assume $l=2k+|S|-2d_S+\varepsilon(k+d_S)-1$.
Theorem \ref{main1} then 
immediately follows from Proposition \ref{calc-KO-degree}
in the cases other than $k+d_S \equiv 2 \pmod 4$,
so we consider this case.
Let $\psi: \HH_1 \to \tilde \R^3$ be
the $\pin$-equivariant map 
defined by $\psi(q)=q i \bar q$, where $\pin$ acts on $\HH_1$
by right multiplication.
Then by considering the direct sum $\varphi \oplus \psi$.
we have the inequality  $l+3 \ge 2(k+1) +|S| -2d_S +3$,
since $(k+1) +d_S \equiv 3 \pmod 4 $.
This completes the proof of Theorem \ref{main1}.

The equation (\ref{Chern})
in Introduction follows immediately
since
$\kappa \sigma \beta_{\HH}(\tilde \R^4)=\beta_{\C}(\C^2)$
by Lemma \ref{forgetful1} and 
\begin{align*}
[V_1]-[W_1]
= k \sigma[\HH] + \sum_{S \subset \{ 1, \cdots, n\}, \caF=4}
a_S \sigma\beta_{\HH}(\tilde \R^S)
+\text{(torsion)}
\in Ksp(\tT^n).
\end{align*}

\section{Proof of Theorem~\ref{main2}}\label{Proof-of-2}

Let $X$ be a connected closed oriented spin $4$-manifold
with indefinite intersection form.
Let $k=-\sign(X)/16$ and $l=b^+(X)>0$.
We write the spinor bundle of $X$
as $S^+(TX) \oplus S^-(TX)=S^+ \oplus S^-$ for simplicity.
Take a Riemannian metric $g$ on $X$.
Then we consider the monopole equation \cite{Witten} as the map:
\begin{align*}
\Tilde\Phi: \Tilde{\cal V}=\Ker (d^*: \Omega^1 \to \Omega^0)\oplus \Gamma(S^+)
& \rightarrow \Tilde{\cal W}=\Omega^+\oplus H^1(X; \R)\oplus \Gamma(S^-),\\
(a, s)& \mapsto P(a,s)+Q(a,s),
\end{align*}
where $P$ is a linear map and $Q$ is a quadratic form
which are constructed as follows:
Let $\pi: \Ker (d^*: \Omega^1 \to \Omega^0) \to 
H^1(X; \R)$ be the projection to the harmonic part
and $D: \Gamma(S^+) \to \Gamma(S^-)$
the Dirac operator of $X$.
Then the maps $P$ and $Q$ are defined to be
$$P(a,s)=(d^+a,\pi(a),Ds), \quad Q(a, s)=(q(s), 0, c(a)s),$$
where $q: S^+ \to \Lambda^+$ is the quadratic map
formed from a nontrivial quadratic $\Spin(4)$-equivariant map
$\Delta_4^+ \to \Lambda^+(\R^4)$.

We should remark that 
if we regard $a$ as a $U(1)$-connection on the trivial complex
line bundle over $X$
then the differential operator $D_a: \Gamma(S^+) \to \Gamma(S^-)$
defined by $D_as=Ds +c(a)s$
is the Dirac operator twisted by $a$,
and $d^+a$ is the self-dual part $F^+_a$ of the curvature
of $a$.

The gauge symmetry of the monopole equation
induces the following symmetry:
Choose a base point $x_0$ in $X$.
Let $\rho: X \to \Hom(H^1(X; \Z) , \R/\Z)$ be
the Albanese map, which is given by
$x \in X \to \{ \omega \to  \int^x_{x_0}\omega \bmod \Z \}$.
Let $\hat \rho: H^1(X; \Z) \to C^{\infty}(X, \R/\Z)$ be 
the adjoint map of $\rho$.
Then the actions of $H^1(X; \Z)$ on 
$\Omega^1, \Gamma(S^{+}\oplus S^-)$, $\Omega^+$ 
are respectively defined by
\begin{align}\label{action}
h\cdot a =a+h, \quad h\cdot s=e^{2\pi\sqrt{-1}\hat \rho(h)}s,
\quad h\cdot b=b  \quad \text{for $h \in H^1(X;\Z)$}.
\end{align}
We introduce the diagonal actions of $H^1(X;\Z)$ 
on $\Tilde{\cal V}$ and $\Tilde{\cal W}$.
Since the projection $\Tilde{\cal V} \to H^1(X;\R)$
is $H^1(X;\Z)$-equivariant,
it descends to the vector bundle ${\cal V}=\Tilde{\cal V}/H^1(X; \Z)$ 
over the Jacobian torus $J_X=H^1(X;\R)/H^1(X;\Z)$
with the decomposition ${\cal V}={\cal V}_0 \oplus {\cal V}_1$ 
into the product real bundle 
${\cal V}_0= J_X \times K$
($K=\Ker d^* \cap H^1(X; \R)^{\perp}$)
and a complex vector bundle ${\cal  V}_1$ 
over $J_X$
with fiber $\Gamma(S^+)$.
Similarly we see
the projection $\Tilde{\cal W} \to H^1(X;\R)$
descends to the vector bundle ${\cal W}=\Tilde{\cal W}/H^1(X;\Z)$
over $J_X$ 
with the decomposition ${\cal W}={\cal W}_0 \oplus {\cal W}_1$
into the product real bundle 
${\cal W}_0= J_X \times \Omega^+$ 
and a complex vector bundle ${\cal  W}_1$ over $J_X$
with fiber $\Gamma(S^-)$.
Then $\Tilde\Phi$ induces a $\dpin$-equivariant
fiber-preserving map
$$\Phi: {\cal V}_0 \oplus {\cal V}_1
\rightarrow {\cal W}_0 \oplus {\cal W}_1.
$$
We should note that the restriction 
$\Phi_0: {\cal V}_0 \to {\cal W}_0$ of $\Phi$
is given by the linear injective map $d^+: K \to \Omega^+$.

The compactness of the moduli space implies
that the inverse image $\Phi^{-1}(0)$ of the zero section
is compact.

Moreover the spin structure of $X$
provides an extra symmetry \cite{Kr}, \cite{Witten}:
The linear involution on $\Omega^1$, $\Omega^+$
and quaternionic scalar multiplication on $\Gamma(S^{\pm})$
induce $\pin$-actions on $\Tilde{\cal V}, \Tilde{\cal W}$
via the projection $\pin \to \{\pm 1\}$ or the inclusion $\pin \to Sp(1)$,
which respectively induce $\pin$-actions on
${\cal V}$, ${\cal W}$, so that $\Phi$ is $\pin$-equivariant.

By a finite-dimensional approximation \cite{Bauer-Furuta}, \cite{Furuta},
we get a proper $\pin$-equivariant fiber-preserving map:
\begin{align*}
\varphi: V_0 \oplus V_1 & \rightarrow W_0 \oplus W_1
\end{align*}
such that 
(i) $V_0=J_X \times \tilde \R^{x}$, $W_0=J_X \times \tilde \R^{x+l}$
for some $x \ge 0 $,
(ii)
$\varphi$ induces the identity on the base space $J_X$,
(iii) the restriction $\varphi_0: V_0 \to W_0$ of $\varphi$
is given by the standard linear inclusion
$\tilde \R^{x} \to \tilde \R^{x+l}$,
(iv) $V_1$, $W_1$ are symplectic bundles over $\tT^n$
whose difference $[V_1]-[W_1] \in Ksp(\tT^n)$
is the index bundle $\Ind \Bbb D$ of the family of Dirac operators
parameterized by $J_X$: 
\begin{align*}
\Bbb D: (H^1(X; \R) \oplus \Gamma(S^+) )/H^1(X; \Z)
& \rightarrow (H^1(X; \R) \oplus \Gamma(S^-) )/H^1(X; \Z),\\
[(a,s)] & \mapsto [(a,D_as)].
\end{align*}

We next calculate $\Ind \Bbb D$ by
identifying it as the index bundle
for a family introduced by G.~Lusztig \cite{Lusztig}:
We take a basis $x_1, \cdots, x_n $ of $H^1(X;\Z)$,
so that 
$J_X$ is identified with the Real torus $\tT^n=(\tilde \R/\Z)^n$.
Let $\hT^n =\R^n /\Z^n$ be $n$-dimensional torus
with the trivial involution.
Let $\Z^n \times \Z^n$ act on $\R^n \times \tilde \R^n \times \C$ by
\begin{align*}
(z_1,z_2)\cdot(r_1,r_2,c)=(r_1+z_1,r_2+z_2,
e^{2\pi \sqrt{-1} \langle z_1, r_2 \rangle}c),
\end{align*}
where $\langle, \rangle$ is the standard inner product on $\R^n$.
The orbit space $L=(\R^n \times \tilde \R^n \times \C)/(\Z^n \times \Z^n)$ is then
a Real line bundle over $T^n \times \tT^n$
with the involution $-1$ defined by
\begin{align*}
(-1)\cdot[(r_1,r_2,c)]=[(r_1, -r_2,\bar c)].
\end{align*}
Put $L_X=(\rho \times \id_{\tT^n})^*L$.
If we denote by $\pi : X \times \tT^n \to \tT^n$ the projection,
we have a $\Z/2$-graded symplectic bundle
$(\pi^*S^{+}\otimes L_X)
\oplus (\pi^*S^{-}\otimes L_X)$
over $X\times \tT^n$
and a symplectic family of Dirac operators
$$ {\Bbb D}' : \Gamma(\pi^*S^{+}\otimes L_X)
\to \Gamma(\pi^*S^{-}\otimes L_X)$$
parameterized by $\tT^n$.
Then it may be obvious from the action (\ref{action})
that $\Ind \Bbb D=\Ind \Bbb D' \in Ksp(\tT^n)$.

Hence, as shown in \cite{Lusztig} and \cite{Ruberman},
the cohomological formula of
the index theorem for families \cite{A-S-IV} indicates
\begin{align*}
\ch(V_1)-\ch(W_1)
& =2k+ \sum_{S \subset \{ 1, \cdots, n\}, \caF=4}
\langle \prod_{i \in S}x_i, [X] \rangle \bigwedge_{i \in S}d\xi_i.
\end{align*}
The equation (\ref{cup-product}) then 
follows from  (\ref{Chern}) in Introduction.

We next describe how the map (\ref{parity}) in Introduction
comes from the index $\Ind \Bbb D$.
Recall that, under the identification $\R^{1,1}=\C$,
$\beta_{\C}(\R^{1,1}) \in KR(\R^{1,1})$ is
the Bott class $\beta_{\C}(\C)$
together with the involution given by complex conjugation \cite{Atiyah-R}.
Then it is easy to show
\begin{align}\label{universal}
[L]
& =\prod_{1 \le q \le n}
\pi_{\{q\}}^*(i_{\{q\}})_!
([\C]+\beta_{\C}(\R^{\{q\},\{q\}})) \notag \\
& =\sum_{S \subset \{1, \cdots, n \}}
\pi_{S}^*(i_{S})_!
\beta_{\C}(\R^{S, S})
\in KR(T^n \times \tT^n),
\end{align}
where $\R^{S, S}=\Map(S, \R^{1,1})$ and
$\beta_{\C}(\R^{S, S}) \in KR(\R^{S, S})$ is the Bott class
corresponding to $\beta_{\C}(\R^{|S|, |S|}) \in KR(\R^{|S|, |S|})$
under the identification $\R^{S, S} \cong \R^{|S|, |S|}$
for each $S \subset \{1, \cdots, n\}$.

Since the symbol class of the Dirac operator $D$
is just the Bott class $\beta_{\HH}(TX) \in KSp(TX)=Ksp(TX)$,
we use the index theorem for families \cite{A-S-V} 
(see also \cite[Section 16]{Lawson}) to get
\begin{align*}
\Ind \Bbb D
=(i \times \id_{\tT^n})_! (\rho \times \id_{\tT^n})^*[L]
\in Ksp (\tT^n) \cong KR^{-4}(\tT^n),
\end{align*}
where $i: X \to \pt$ is the constant map.
Substituting (\ref{universal}) for $[L]$
and applying the periodicity isomorphisms 
$$KR(\R^{S,S})\cong KO^{|S|}(\R^S),
\quad KR^{-4}(\tilde \R^S)\cong KO^{|S|+4}(\pt)$$
(see Section \ref{real-symplectic}), 
we deduce that $a_S$ ($|S|=2$ or $3$) is zero if and only if
the map (\ref{parity}) is zero.

  {\scriptsize
  Graduate School of Mathematical Sciences,
  University of Tokyo, Tokyo 153-8914, Japan
  (furuta@ms.u-tokyo.ac.jp)}

  {\scriptsize Department of Mathematics,
  Keio University, Yokohama 223-8522, Japan
  (kametani@math.keio.ac.jp)}


\begin{thebibliography}{11}

\bibitem{Adams} J.~F.~Adams, 
{\em Lectures on Lie group},
W. A. Benjamin, Inc., (1969).

\bibitem{Atiyah-R} M.~F.~Atiyah, 
{\em $K$-theory and Reality},
Quart. J. Math., Oxford(2), {\bf 17} (1966), 367-386.

\bibitem{Atiyah-K} M.~F.~Atiyah, 
{\em $K$-theory},
Benjamin, New York (1967).

\bibitem{Atiyah} M.~F.~Atiyah, 
{\em Bott periodicity and the index of elliptic operators},
Quart. J. Math., Oxford(2), {\bf 19} (1968), 113-140.

\bibitem{Atiyah-B-S} M.~F.~Atiyah, R.~Bott and A.~Shapiro,
{\em Clifford modules},
Topology, {\bf 3} (suppl. 1) (1964), 3-38.

\bibitem{A-S-IV} M.~F.~Atiyah and I.~M.~Singer, 
{\em The index of elliptic operators IV},
Ann. of Math., {\bf 93} (1971), 119-138.

\bibitem{A-S-V} M.~F.~Atiyah and I.~M.~Singer, 
{\em The index of elliptic operators V},
Ann. of Math., {\bf 93} (1971), 139-149.

\bibitem{AtiyahTall} M.~F.~Atiyah and D.~O.~Tall,
{\em Group representations, $\lambda$-ring and the $J$-homomorphism},
Topology, {\bf 8} (1969), 253-297.


\bibitem{Bauer-Furuta} S.~Bauer and M.~Furuta,
{\em A stable cohomotopy refinement of Seiberg-Witten invariansts, I},
Invent. Math, {\bf 155} (2004), 1-19.

  \bibitem{Donovan-Karoubi} P.~Donovan and M.~Karoubi,
{\em Graded Braner groups and $K$-theory with local coefficients}, 
Publ. Inst. Hautes \'{E}tudes Sci. Publ. Math., No. 38
(1970), 5-25.

\bibitem{Dupont} J.~L.~Dupont,
{\em Symplectic bundles and $KR$-theory},
Math. Scand., {\bf 24} (1969), 27-30.

  \bibitem{Furuta} M.~Furuta, 
{\em Monopole equation  and  the $11/8$ conjecture},
Math. Res. Lett., {\bf 8} (2001), 279-291.

  \bibitem{FKMa} M.~Furuta, Y.~Kametani and H.~Matsue,
  {\em Spin 4-manifolds with signature=-32},
Math. Res. Lett., {\bf 8} (2001), 293-301.

  \bibitem{FKMM} M.~Furuta, Y.~Kametani,
  H.~Matsue and N.~Minami,
  {\em Stable-homotopy Seiberg-Witten invariants
  and Pin bordisms},
  preprint (UTMS 2000-46, University of Tokyo).

  \bibitem{FKM} M.~Furuta, Y.~Kametani and N.~Minami,
  {\em Stable-homotopy Seiberg-Witten invariants
  for rational cohomology $K3 \# K3$'s},
  J. Math. Sci. Univ. Tokyo, {\bf 8} (2001), 157-176.

  \bibitem{Husemoller} D.~Husemoller,
  {\em Fiber bundles, Third Edition},
  Graduate Texts in Mathematics {\bf 20} Springer 
Verlag, New York, (1994).

\bibitem{Komiya} K.~Komiya,
{\em Equivariant $K$-theory and maps between representation spheres},
Publ. Res. Inst. Math. Sci., {\bf 31} (1995), 725-730.

  \bibitem{Kr} P.~B.~Kronheimer,
lecture at Cambridge University  in December 1994.

\bibitem{Lawson} H.~B.~Lawson and M.~Michelson, 
{\em Spin Geometry}, Princeton Univ. Press, (1989).

  \bibitem{Lusztig} G.~Lusztig, {\em Novikov's higher signature
and families of elliptic operators}, J. Differential Geom.,
{\bf 7} (1972), 229--256.

  \bibitem{Minami} N.~Minami, {\em The $G$-join theorem
  - an unbased $G$-Freudenthal theorem}, preprint.

\bibitem{Ruberman} D.~Ruberman and S.~Strle, {\em \ Mod 2 Seiberg-Witten
invariants of homology tori},
Math. Res. Lett., {\bf 7} (2000), 789-799.


 \bibitem{Schmit} B.~Schmidt, {\em Spin 4-manifolds and {\rm Pin(2)}-equivariant
homotopy theory}, Ph.D thesis, Universit\"{a}t Bielefeld, (2003).

\bibitem{Segal} G.~B.~Segal, {\em Equivariant K-theory},
Publ. Inst. Hautes \'{E}tudes Sci. Publ. Math., No. 34
(1968), 129-151.


  \bibitem{Stolz} S.~Stolz,
  {\em The level of real projective spaces},
  Comment. Math. Helvetici {\bf 64} (1989), 661-674.

  \bibitem{Witten} E.~Witten, 
{\em Monopoles and four-manifolds},
Math. Res. Lett., {\bf 1} (1994), 769-796.


  \end{thebibliography}
\end{document}